\documentclass[10.5pt]{article}

\usepackage{amsmath,amsfonts,amsbsy,amstext,amsopn,amssymb,bm,graphicx,epstopdf,float}
\usepackage{longtable}
\usepackage{theorem}
\usepackage{mathrsfs}
\usepackage{epsfig}
\usepackage{indentfirst}
\usepackage[noend]{algpseudocode}
\usepackage{algorithmicx,algorithm}
\usepackage{algpseudocode}
\usepackage{lipsum}
\usepackage{subfig}
\usepackage{verbatim}
\usepackage{cite}
\usepackage{natbib}
\usepackage{multirow}
\makeatletter \@addtoreset{equation}{section}

\makeatother

{\theorembodyfont{\slshape}

} {\theorembodyfont{\rmfamily}

}

\usepackage{graphicx}

\usepackage{color}
\setcitestyle{numbers,square}
\begin{document}
\title{\bf An Augmented Subspace Based Adaptive Proper Orthogonal Decomposition Method for Time Dependent Partial Differential Equations}
\author{Xiaoying Dai\thanks{LSEC, Institute of Computational Mathematics and Scientific/Engineering Computing,
Academy of Mathematics and Systems Science, Chinese Academy of Sciences,  Beijing 100190, China; and School of Mathematical Sciences,
University of Chinese Academy of Sciences, Beijing 100049, China
{E-mail: {\tt daixy@lsec.cc.ac.cn}}}
\and Miao Hu\thanks{LSEC, Institute of Computational Mathematics and Scientific/Engineering Computing,
Academy of Mathematics and Systems Science, Chinese Academy of Sciences,  Beijing 100190, China; and School of Mathematical Sciences,
University of Chinese Academy of Sciences, Beijing 100049, China
{E-mail: {\tt humiao@lsec.cc.ac.cn}}}
\and Jack Xin\thanks{Department of Mathematics, University of California at Irvine, Irvine, CA 92697, USA
{E-mail: {\tt jack.xin@uci.edu}}}
\and Aihui Zhou\thanks{LSEC, Institute of Computational Mathematics and Scientific/Engineering Computing,
Academy of Mathematics and Systems Science, Chinese Academy of Sciences,  Beijing 100190, China; and School of Mathematical Sciences,
University of Chinese Academy of Sciences, Beijing 100049, China
{E-mail: {\tt azhou@lsec.cc.ac.cn}}}
}
\date{}
\maketitle

\begin{quote}
    {\small {\bf Abstract}\\

        In this paper, we propose an augmented subspace based adaptive proper orthogonal decomposition (POD) method for solving the time dependent partial differential equations. By augmenting the POD subspace with some auxiliary modes, we obtain an augmented subspace. We use the difference between the approximation obtained in this augmented subspace and that obtained in the original POD subspace to construct an error indicator, by which we obtain a general framework for augmented subspace based adaptive POD method. We then provide two strategies to obtain some specific augmented subspaces, the random vector based augmented subspace and the coarse-grid approximations based augmented subspace. We apply our new method to two typical 3D advection-diffusion equations with the advection being the Kolmogorov flow and the ABC flow. Numerical results show that our method is more efficient than the existing adaptive POD methods, especially for the advection dominated models.
    }


    \date{Received: date / Accepted: date}
\end{quote}

{\small {\bf Keywords} \\

Proper orthogonal decomposition $\cdot$ Adaptive $\cdot$ Augmented subspace $\cdot$  Galerkin projection $\cdot$ Error indicator
}

\section{Introduction}
Time dependent partial differential equations arise in many important fields, e.g., the seawater intrusion \cite{Bakker2000}, the heat transfer \cite{Cannon1984}, semiconductor devices \cite{Markowych1990} and fluid equations \cite{Holmes1996, Batchelor1999, Xin2005, Zu2015, Xin2021flow, Xin2022flow}. The study about the numerical methods for the time dependent partial differential equations is an important and attractive research topic. There are some classical numerical discretization methods for the spatial discretization of the time dependent partial differential equations, e.g., the finite element method \cite{Brenner2008}, the finite difference method \cite{Strikwerda2004} and the plane wave method \cite{Kosloff1983}. Usually, the semi-discretized systems resulted from applying these classical spatial  discertization  methods for a time dependent partial differential equations are huge dimensional. If we use these classical spatial discretization methods at each time interval, the computational cost will be very expensive, especially for complex systems.

Therefore, some efficient and accurate model order reduction methods have been proposed to reduce the dimension of discretized system and then the computational costs \cite{Benner2015, Quarteroni2014, Chinesta2017, maday2006reduced, boyaval2010reduced}.  The basic idea for the model order reduction method is to project the continuous system onto a low-dimensional approximation subspace whose dimension is significantly less than that of the classical discretization space. Proper orthogonal decomposition (POD) is a commonly used model order reduction technique \cite{Hesthaven2016, Volkwein2010, Pinnau2008, APODXin2017, Grable2018}.  The typical steps for getting a POD reduced order model of a time dependent partial differential equation are as follows. First, choose  some classical discretized method to discretize the continuous equation for some time interval, and solve the resulted high-dimensional  discretized systems to get a set of snapshots. Then, construct the POD modes by minimizing the error between these modes and the snapshots, which is equivalent to solving an $N_s$ dimensional eigenvalue problem  \cite{Sirovich1987, Holmes1996,Rowley2004}. Here $N_s$ is the number of snapshots. At last, project the original system on the subspace spanned by the POD modes, which is also called POD subspace, and solve the resulted discretized systems  in the following time intervals. In the actual calculation,  singular value decomposition (SVD) is usually used to obtain the POD modes from the snapshots \cite{Ly2001, Pinnau2008,burkardt2006pod}. By choosing the snapshots properly, the number of the POD modes will usually be of magnitude smaller than the degrees of freedom resulted by the classical spacial discretization  methods.

 There are many applications of the POD methods in scientific and engineering field. For instance, in \cite{Dickinson2010}, a group proper orthogonal decomposition (GPOD) method was introduced to simulate the nonlinear Burgers equation. In \cite{Li2019}, the POD method was used to solve the time-domain Maxwell's equations. In \cite{Li2021}, a reduced-order finite element formulation based on POD method was established for the Allen-Cahn equation. Other applications include studies of turbulence \cite{Druault2005,Tabib2008}, process identification \cite{Katayama2005, Khalil2007} and control in chemical engineering \cite{Padhi2003,ChaoXu2011}, etc. We refer to \cite{Volkwein2011, Volkwein2013} for more introduction to the POD method.

There are also some existing works on the error analysis for the POD method. For example, Kunisch and Volkwein estimated the error of the POD approximation for linear and nonlinear evolution equations in \cite{Volkwein2001,Volkwein2002}. In \cite{Xin2021}, Xin et al. analyzed the convergence of the POD approximation for viscous G-equation by decomposing the data into a mean-free part and a mean part. More works about the numerical analysis for the POD method can be referred to \cite{Homescu2005, Locke2020,Koc2021} and references therein. 

The classical POD method for time dependent partial differential equations only uses the snapshots obtained in the early time interval to construct the POD modes. Once the POD modes are obtained, they will be fixed and not updated during the time evolution any more. However, the solution of the system may change a lot over time. Therefore, if the POD modes are not updated as the time evolution, the approximation error obtained by the POD methods may  become larger and larger.

In order to improve the accuracy of the POD methods in the whole time interval, some adaptive POD methods which update the POD modes as time evolution have been introduced in recent years \cite{Dai2017,Kuang2020, Rapun2015, Terragni2015, Peherstorfer2015}. In \cite{Dai2017, Rapun2015, Terragni2015}, the authors constructed some residual type error estimators, based on which some residual based adaptive POD methods are proposed for simulation of time dependent problems. In \cite{Kuang2020}, the authors proposed a two-grid based adaptive POD method. For this method, they first constructed two finite element spaces, a coarse finite element space and a fine finite element space, and then used the error obtained in the coarse finite element space to construct the error indicator, by which people then knew if it is needed to update the POD subspace in the fine mesh or not. We refer to \cite{Kuang2020} for more details about the two-grid based adaptive POD method.

In this paper, we propose a new approach for developing some adaptive POD methods for the time dependent partial differential equations. The main idea of our approach is to add some auxiliary modes to the POD subspace to build an augmented subspace. We then use the difference between the approximation obtained in the augmented subspace and that obtained in original POD subspace to develop an error indicator. Using this idea, we obtain a general framework for augmented subspace based adaptive POD method. We then provide two specific methods to obtain the auxiliary mode, one is simply using a randomly generated vector as the auxiliary mode, the other is using a coarse-grid approximation as the auxiliary mode. Using the coarse-grid approximation to build the augmented subspace is inspired by \cite{Kuang2020}. In \cite{Kuang2020}, the authors used the error obtained in a coarse finite element space to act as the error indicator, while in this paper, we use the approximation in the coarse finite element space to augment the POD subspace, and then develop an error indicator.

The rest of this paper is organized as follows. In Section 2, we give some preliminaries, including the general framework of the adaptive POD method, the complexity of the POD type methods and the simple descriptions of two typical adaptive POD methods. In Section 3, we propose a general framework for the  augmented subspace based adaptive POD method, and provide two specific methods to construct the augmented subspace, the random vector based augmented subspace and the coarse grid approximations based augmented subspace.  In Section 4, we apply the coarse-grid approximations based adaptive POD method  to the simulation of some typical time dependent partial differential equations, i.e., the advection-diffusion equation with three dimensional velocity field, including both the Kolmogorov flow and the ABC flow. The numerical results show the accuracy and efficiency of our new method. In Section 5, we give some concluding remarks. Finally, we provide some additional numerical results for the advection dominated models with different coarse meshes and different error indicator threshold in ``Appendix A''.

\section{Preliminaries}
 We first recall some definition and notation. We shall use the standard notation for Sobolev spaces  and their associated norms and seminorms; see, e.g., \cite{Adams1975}. Let $V$ be a Banach space with norm $\Vert*\Vert_V$ and $L^{p}(0, T; V)$ be a Banach space equipped with the norm 
\begin{equation*}
    \Vert u \Vert_{ L^{p}(0, T; V)} = \left(\int^T_0 \Vert u(t)\Vert_{V}^p \textup{d}t\right)^\frac{1}{p}, 1\leq p <\infty.
\end{equation*}
We call $u\in  C(0, T; V)$ if
\begin{equation*}
    \Vert u \Vert_{C(0, T; V)}=\mathop{\max}\limits_{t\in[0, T]} \Vert u(t)\Vert_{V}<\infty.
\end{equation*}
In this paper, vectors and matrices will be denoted by bold letters.

We consider the following general time dependent partial differential equations:
\begin{align}\label{equ}
    \left\{
    \begin{aligned}
         & u_{t} - \epsilon\Delta u + {\bf B}(x, y, z, t)\cdot \nabla u + c(x, y, z, t)u = f(x, y, z, t), \quad \text{in} ~\Omega\times (0, T] \\
         & u(x, y, z, 0) = h(x, y, z),                                                                                                       \\
         & u(x+\kappa, y, z, t) = u(x, y+\kappa, z, t) = u(x, y, z+\kappa, t) = u(x, y, z, t),
    \end{aligned}
    \right.
\end{align}
where $\Omega=[0, \kappa]^{3}$, $f\in L^{2}(0, T; L^{2}(\Omega))$, $c\in L^{\infty}(\Omega)$, ${\bf B}\in C(0, T; W^{1, \infty}(\Omega)^{3})$ and $\epsilon$ is a constant.
We take an inner product with $v$ in (\ref{equ}) and define a bilinear form to simplify the variational form
\begin{equation*}
    a(t; u, v) = \epsilon(\nabla u, \nabla v) - \epsilon\int_{\partial\Omega} \frac{\partial u}{\partial n} vd\sigma + ({\bf B}\cdot \nabla u, v) + (cu, v), \forall u, v\in H^{1}(\Omega),
\end{equation*}
where $(\cdot, \cdot)$ stands for the inner product in $L^{2}(\Omega)$.
We can obtain the variational form of the Eq. (\ref{equ}) as follows: find $u\in L^{2}(0, T; V)$, $u_{t}\in L^{2}(0, T; V^{*})$ such that
\begin{align}\label{weak}
    (\frac{\partial u}{\partial t}, v) + a(t; u, v) = (f(x, y, z, t), v), \forall v\in V,
\end{align}
where
$$V = \{v\in H^{1}(\Omega) : v|_{x=0} = v|_{x=\kappa}, v|_{y=0} = v|_{y=\kappa}, v|_{z=0} = v|_{z=\kappa}\},$$
and $V^{*}$ is the dual space of V.
\subsection{Standard discretization}
We first consider the standard temporal discretization of (\ref{weak}). There are many existing temporal discretization methods, such as Euler method and implicit Euler method \cite{He2008}, which can be used to discretize (\ref{weak}). Here, we choose the implicit Euler method. We first divide the time interval into $N\in \mathbb{N}$ subintervals with equal length $\delta t = T/N$, and set $u^{k}(x, y, z) = u(x, y, z, t_{k})$, $f_{k}(x, y, z) = f(x, y, z, t_{k})$, where $t_{k} = k*\delta t$, for $k\in \{0, 1, \cdots, N\}.$ Then we can get the semi-discretization scheme of (\ref{weak}) as follows:
\begin{align}\label{semi}
    \left(\frac{u^{k}(x, y, z)-u^{k-1}(x, y, z)}{\delta t}, v\right) + a(t_{k}; u^{k}(x, y, z), v) = (f_{k}(x, y ,z), v), \quad v\in V.
\end{align}

We then consider the classical spatial discretization of (\ref{semi}). Here, we choose the finite element method to discretize (\ref{semi}). Let $\mathcal{T}_h$ be a shape regular family of nested
conforming mesh over $\Omega$ with size $h$:
there exists a constant $\gamma^{\ast}$ such that
\begin{eqnarray}
    \frac{h_{\tau}}{\rho_{\tau}} \leq \gamma^{\ast}, \quad\forall~\tau \in \mathcal{T}_h,
\end{eqnarray}
where $h_{\tau}$ is the diameter of $\tau$ for each $\tau\in \mathcal{T}_h$, $\rho_{\tau}$ is the
diameter of the biggest ball contained in $\tau$,
and $h=\max\{h_{\tau}: \tau\in \mathcal{T}_h\}$. Let $V_h$ be a subspace of continuous functions on $\Omega$ such that
\begin{align*}
    V_{h} = \{v_{h} : v_{h}|_{\tau} \in \mathbb{P}_{\tau}, \forall~\tau\in \mathcal{T}_{h}  \mbox{ and }  v_{h} \in C^{0}(\overline{\Omega})\},
\end{align*}
where $\mathbb{P}_{\tau}$ is a set of polynomials on element $\tau$. Let $\{\phi_{h,1}, \phi_{h,2},\cdots,\phi_{h,N_g}\}$ be a basis of $V_h$. Denote
\begin{align*}
    \Phi_{h}:= (\phi_{h,1}, \phi_{h,2},\cdots,\phi_{h,N_g}).
\end{align*}
The standard finite element discretization of (\ref{semi}) can be formulated as follows:
find $u^{k}_{h}(x, y, z)\in V_h$, such that
\begin{align}\label{femsemi}
    \left(\frac{u^{k}_h(x, y, z)-u^{k-1}_h(x, y, z)}{\delta t}, v_h\right) + a(t_{k}; u_h^{k}(x, y, z), v_h) = (f_{k}(x, y ,z), v_h), \quad v_h\in V_h.
\end{align}
Note that $u_h^k(x,y,z)$ can be expressed as
\begin{eqnarray}\label{discform}
    u^{k}_{h}(x, y, z)=\mathop{\sum}^{N_g}\limits_{i=1} u_{h,i}^{k}\phi_{h,i}(x, y, z).
\end{eqnarray}
Inserting (\ref{discform}) into (\ref{femsemi}), and setting $v_h=\phi_{h,j}, j=1, 2, \cdots,N_g$, respectively, we have
\begin{eqnarray}\label{femfull}
    \left( \sum^{N_g}_{i=1} u_{h,i}^{k} \phi_{h,i} -\sum^{N_g}_{i=1}u_{h,i}^{k-1} \phi_{h,i}, \phi_{h, j}\right) + \delta t a(t_{k}; \sum^{N_g}_{i=1} u_{h,i}^{k} \phi_{h,i}, \phi_{h, j}) = \delta t(f_{k}, \phi_{h, j}).
\end{eqnarray}
Define
\begin{align*}
    {\bf A}^{k}_{h, ij} = (\phi_{h,j}, \phi_{h,i}) +\delta t a(t_{k}; \phi_{h,j}, \phi_{h,i}), \quad & {\bf u}^{k}_{h}=(u_{h,1}^{k}, u_{h,2}^{k}, \cdots, u_{h,N_g}^{k})^{T}, \\
    {\bf b}^{k}_{h}=\delta t * ((f_k, \phi_{h,1}),\cdots,(f_k,\phi_{h,N_g}))^{T}, \quad              & {\bf C}_{h, ij}=(\phi_{h,j}, \phi_{h,i}).
\end{align*}
Then (\ref{femfull}) can be rewritten as the following algebraic form
\begin{align}\label{FEMform}
    {\bf A}^{k}_{h} {\bf u}^{k}_{h} = {\bf b}^{k}_{h} +{\bf C}_{h} {\bf u}^{k-1}_{h}.
\end{align}
\subsection{Adaptive POD method}
In this subsection, we first recall the general procedure for getting a POD reduced order model of a time dependent partial differential equation. First, we choose some classical discretized method to discretize (\ref{semi}) in some time interval. Here we discretize (\ref{semi}) in the finite element space  $V_{h}$ for $t\in[0, T_{0}]$, and solve the resulted high dimensional discretized systems at different times $t_{0}$, $t_{\delta M}$, $\cdots$, $t_{n_{s}\cdot \delta M}$ to get the snapshot matrix ${\bf U}_{h}$. Here, $\delta M$ is an integer parameter and $n_{s}=\lfloor\frac{T_{0}}{\delta t\cdot\delta M}\rfloor$, where $\lfloor * \rfloor$ means the round down. Then, we perform SVD on ${\bf U}_{h}$, and obtain ${\bf U}_{h} = {\bf RSV^{T}}$. Note that the diagonal elements in {\bf S} are arranged from largest to smallest, we set the number of POD modes by
\begin{align}\label{PODnum}
    m = \min\{k| \mathop{\sum}^{k}\limits_{i=1} {\bf S}_{i, i} > \gamma_{1}* {\rm Trace}({\bf S})\},
\end{align}
where $\gamma_{1}$ is a given parameter, and set  ${\bf\widetilde R} = {\bf R}[:, 1:m]$. Then the POD modes are
$$\Psi_{h} = (\psi_{h,1}, \psi_{h,2}, \cdots, \psi_{h, m}) := \Phi_{h}{\bf\widetilde R}.$$
The process for constructing POD modes can be summarized as a routine POD\_Mode(${\bf U}_{h}, \gamma_1, \Phi_{h}, m, \Psi_{h}$) in Algorithm \ref{PODMode},  see also \cite{Kuang2020}.
\begin{algorithm}[t]
    \caption{POD\_Mode(${\bf U}_{h}, \gamma, \Phi_{h}, m, \Psi_{h}$)}\label{PODMode}
    \hspace*{0.02in} {\bf Input:} ${\bf U}_{h}, \gamma, \Phi_{h} = (\phi_{h,1}, \phi_{h,2}, \cdots, \phi_{h,n})$.\\
    \hspace*{0.02in} {\bf Output:} m and POD modes $\Psi_h = (\psi_{h,1}, \psi_{h,2}, \cdots, \psi_{h, m})$.
    \begin{algorithmic}[1]
        \State Perform SVD on ${\bf U}_{h}$ to obtain ${\bf U}_{h} = {\bf RSV^{T}}$, where {\bf S} = diag $\{\sigma_{1}, \sigma_{2}, \cdots, \sigma_{r} \}$ with $\sigma_{1} \geq \sigma_{2}\geq\cdots\geq\sigma_{r} > 0.$
        \State Set $m = min\{k| \mathop{\sum}^{k}\limits_{i=1} {\bf S}_{i,i} > \gamma* {\rm Trace}({\bf S})\}.$
        \State $(\psi_{h,1}, \psi_{h,2}, \cdots, \psi_{h, m}) = \Phi_{h}{\bf R}[:, 1:m]$.
    \end{algorithmic}
\end{algorithm}
At last, we project (\ref{semi}) on the POD subspace $V_{h, \rm{POD}} = \text{span}\{\psi_{h,1}, \cdots,\psi_{h,m}\}$ when $t > T_{0}$. The POD approximation  $u^{k}_{h, \rm{POD}}$ is then expressed as
\begin{align}\label{uhPOD}
    u^{k}_{h, \rm{POD}}(x, y, z) = \mathop{\sum}^{m}\limits_{i=1} \alpha_{h,i}^{k}\psi_{h,i}(x, y, z).
\end{align}
Inserting (\ref{uhPOD}) into (\ref{semi}), and setting $v=\psi_{h,j}, j=1, 2, \cdots, m$, respectively, we obtain
\begin{align}\label{podfull}
    \left( \sum^{m}_{i=1} \alpha_{h,i}^{k} \psi_{h,i} -\sum^{m}_{i=1}\alpha_{h,i}^{k-1} \psi_{h,i}, \psi_{h, j}\right) + \delta t a(t_{k}; \sum^{m}_{i=1} \alpha_{h,i}^{k} \psi_{h,i}, \psi_{h, j}) = \delta t(f_{k}, \psi_{h, j}).
\end{align}
Define
\begin{align*}
    {\bf \bar A}^{k}_{h, ij} = (\psi_{h, j}, \psi_{h,i}) +\delta t a(t_{k}; \psi_{h, j}, \psi_{h,i}), \quad & {\bf u}^{k}_{h, \rm{POD}}=(\alpha_{h,1}^{k}, \alpha_{h,2}^{k},\cdots, \alpha_{h,m}^{k})^{T}, \\
    {\bf \bar b}^{k}_{h}=\delta t * ((f_{k}, \psi_{h,1}),\cdots,(f_{k},\psi_{h,m}))^{T}, \quad              & {\bf \bar C}_{h, ij}=(\psi_{h,j}, \psi_{h,i}).
\end{align*}
Then (\ref{podfull}) can be rewritten as the following algebraic form
\begin{align}\label{PODform}
    {\bf \bar A}^{k}_{h} {\bf u}^{k}_{h, \rm{POD}} = {\bf \bar b}^{k}_{h} +{\bf \bar C}_{h} {\bf u}^{k-1}_{h, \rm{POD}}.
\end{align}
From the expression of the POD modes, i.e., $\Psi_{h} =\Phi_{h}{\bf\widetilde R}$, the above equation can also be written as
$$
    {\bf \widetilde{R}}^{T}{\bf A}^{k}_{h}{\bf \widetilde{R}} {\bf{u}}^{k}_{h, \rm{POD}}
        ={\bf \widetilde{R}}^{T}{\bf b}^{k}_{h}
    +
    {\bf \widetilde{R}}^{T}{\bf C}_{h}{\bf \widetilde{R}}
    {\bf u}^{k-1}_{h, \rm{POD}}.
$$

As we mentioned before, the classical POD method only uses the snapshots obtained in time interval $[0,T_0]$ to construct the POD modes. Once the POD modes are obtained, they will not be updated during the time evolution. However, the solution may change a lot over time. In order to improve the accuracy of the POD method in the whole time interval, more and more researchers have paid more attention to the study of adaptive POD method in recent years \cite{Dai2017,Kuang2020, Rapun2015, Terragni2015, Peherstorfer2015}. Here, we provide a briefly introduction to it. 

Motivated by the procedure for adaptive finite element method \cite{Dai2008}, in \cite{Kuang2020}, the authors summarized the procedure of adaptive POD method as a loop constructed by the following four steps:
\begin{enumerate}
    \item \textbf{Solve}: Solve the Eq. (\ref{semi}) in the POD subspace $\text{span}\{\psi_{h,1}, \cdots,\psi_{h,m}\}$.
    \item \textbf{Estimate}: Construct an error indicator $\eta_k$ to estimate the error of the approximation obtained in current POD subspace.
    \item \textbf{Mark}: Mark the time instance $t_k$ when the POD subspace is needed to be updated.
    \item \textbf{Update}: Update the POD subspace at the marked time instance.
\end{enumerate}

The step \textbf{Solve} is just the procedure for obtaining the approximations by the classical POD method we introduced above. The step \textbf{Estimate} is crucial for an adaptive POD method, which determines the efficiency and accuracy of the method. The step \textbf{Mark} picks out the time instance when the POD modes are needed to be updated. For the step \textbf{Update}, it is worthy of noting that, if the time instance $t = q\delta t$ is marked, we will go back to the previous time instance $t_{1} = (q-1)\delta t$ to restart the collection of the approximations in the finite element space and get the snapshots matrix ${\bf W}_{h,1}$. In order to obtain the new POD modes, we perform SVD on ${\bf W}_{h,1}$ to obtain  ${\bf W}_{h,1} = {\bf R}_{1} {\bf S}_{1} {\bf V}_{1}^{T}$. Then, we obtain the number of POD modes $m_1$ by (\ref{PODnum}) but with a different parameter $\gamma_2$. In order to keep the information in original POD modes, we perform SVD on ${\bf W}_{h,2} = [{\bf R}_{1}[:, 1 : m_{1}], {\bf \widetilde R}]$, and get ${\bf W}_{h,2} = {\bf R}_{2} {\bf S}_{2} {\bf V}_{2}^{T}$. Then, we get the number of POD modes $m$ by (\ref{PODnum}) but with a parameter $\gamma_3$. Finally we set ${\bf \widetilde R}={\bf R}_{2}[:, 1 : m]$ and obtain the update POD modes by
\begin{align*}
    \Psi_{h} = (\psi_{h,1}, \psi_{h,2}, \cdots, \psi_{h,m}) = \Phi_{h} {\bf \widetilde R}.
\end{align*}
For the convenience of the following discussion, we summarize the process for the step \textbf{Update} as routine Update\_POD\_Mode(${\bf W}_{h, 1}, \gamma_{2}, \gamma_{3}, \Phi_{h},  m, \Psi_{h}$) in Algorithm \ref{UpdatePODAlg}, which is first introduced in \cite{Kuang2020}.

\begin{algorithm}[H]
    \caption{Update\_POD\_Mode(${\bf W}_{h, 1}, \gamma_{2}, \gamma_{3}, \Phi_{h}, m, \Psi_{h}$)}\label{UpdatePODAlg}
    {\bf Input:} ${\bf W}_{h, 1}, \gamma_{2}, \gamma_{3},  \Phi_{h} = (\phi_{h, 1}, \phi_{h, 2}, \cdots, \phi_{h, n})$, $m$ and $m$ old POD modes $\Psi_{h}, \Psi_{h} = \Phi_{h}{\bf \widetilde R}$.\\
    {\bf Output:} new $m$ and new $m$ POD modes $\{\psi_{h, 1}, \psi_{h, 2}, \cdots, \psi_{h, m}$\}.
    \begin{algorithmic}[1]
        \State Perform SVD on ${\bf W}_{h,1}$ to obtain  ${\bf W}_{h,1} = {\bf R}_{1} {\bf S}_{1} {\bf V}_{1}^{T}$, where ${\bf S_{1}}$ = diag$\{\sigma_{1,1}, \sigma_{1,2}, \cdots, \sigma_{1,r_{1}} \}$ with $\sigma_{1,1} \geq \sigma_{1,2}\geq\cdots\geq\sigma_{1, r_{1}} > 0.$
        \State Set $m_{1} = {\rm min}\{k| \mathop{\sum}^{k}\limits_{i=1} {\bf S}_{1,ii} > \gamma_{2} * {\rm Trace}({\bf S}_{1})\}.$
        \State Perform SVD on ${\bf W}_{h,2} = [{\bf R}_{1}[:, 1 : m_{1}], {\bf \widetilde R}]$, and obtain  ${\bf W}_{h,2} = {\bf R}_{2} {\bf S}_{2} {\bf V}_{2}^{T}$, where ${\bf S_{2}}$ = diag$\{\sigma_{2,1}, \sigma_{2,2}, \cdots, \sigma_{2,r_{2}} \}$ with $\sigma_{2,1} \geq \sigma_{2,2}\geq\cdots\geq\sigma_{2, r_{2}} > 0.$
        \State Set $m = {\rm min}\{k| \mathop{\sum}^{k}\limits_{i=1} {\bf S}_{2,ii} > \gamma_{3} * {\rm Trace}({\bf S}_{2})\}$, and ${\bf \widetilde R} = {\bf R}_{2}[:, 1 : m]$.
        \State $(\psi_{h,1}, \psi_{h,2}, \cdots, \psi_{h,m}) = \Phi_{h}{\bf \widetilde R}.$

    \end{algorithmic}
\end{algorithm}
We then obtain the framework of the adaptive POD method as Algorithm \ref{APODAlg}.
\begin{algorithm}[H]
    \caption{Framework of the adaptive POD method}\label{APODAlg}
    \begin{algorithmic}[1]
        \State Given $T_{0}$, $\delta M$, $\delta T$, $\delta t$,$\gamma _{1}$, $\gamma_{2}$, $\gamma _{3}$, $\eta_0$ and the mesh $\mathcal{T}_{h}$.
        \State Discretize (\ref{semi}) in the standard finite element space $V_h$ on interval $[0, T_{0}]$ and obtain the snapshots matrix ${\bf U}_{h}$.
        \State Construct POD modes $\Psi_{h}$ by POD\_Mode(${\bf U}_{h}, \gamma_1, \Phi_{h}, m, \Psi_{h}$).
        \State $t = T_{0}$, $k = \frac{T_{0}}{\delta t}$.
        \While {$t \leq T$}
        \State  $t = t + \delta t$, k = k + 1.
        \State  Discretize (\ref{semi}) in the POD subspace $V_{h, \rm{POD}}=\text{span}\{\psi_{h,1}, \psi_{h,2}, \cdots, \psi_{h,m}\}$, and obtain the POD approximations ${\bf u}^{k}_{h, \rm{POD}}$.
        \State Compute the error indicator $\eta_k$ by some strategy.
        \If {$\eta_k>\eta_0$}
        \State  $t = t - \delta t$, k = k - 1.
        \State Discretize (\ref{semi}) in $V_{h}$ on interval $[t, t+\delta T]$ and get snapshots ${\bf W}_{h, 1}$, then update POD modes $\Psi_{h}$  by Update\_POD\_Mode(${\bf W}_{h, 1}, \gamma_{2}, \gamma_{3}, \Phi_{h}, m, \Psi_{h}$).
        \State $k = k + \frac{\delta T}{\delta t}$.
        \EndIf
        \State {\bf end if}
        \EndWhile
        \State{\bf end while}
    \end{algorithmic}
\end{algorithm}
\subsection{Complexity}
In this subsection, we analyze the computational complexity of the standard finite element method, the standard POD method and the adaptive POD method introduced before. By comparing the computational complexity of these methods, the advantages of the POD type methods will be shown. We use $\mathcal{O}(n)$ to represent the complexity of a function that increases linearly with respect to $n$. For simplicity, we only consider the case that ${\bf B}$, $c$ and $f$ are separable in time and space, that is, ${\bf B}$ can be denoted as ${\bf B}(x, y, z, t)= {\bf B_1}(x, y, z) + {\bf B_2}(x, y, z) B_3(t)$.

For the standard finite element method, the main computational costs at each time instance are those for building the discretized system and solving the discretized system. We first consider the cost for building the discretized system. We note that during the building of the linear system (\ref{FEMform}), only the terms $({\bf B} \cdot \nabla \phi_{h, j}, \phi_{h,i})$, $(c\phi_{h, j}, \phi_{h,i})$ and $(f, \phi_{h,i})$ are changed as the evolution of time. Let $N_g$ denote the degree of freedom for the finite element discretization. We only need to compute $({\bf B_1} \cdot \nabla \phi_{h, j}, \phi_{h,i})$ and $({\bf B_2} \cdot \nabla \phi_{h, j}, \phi_{h,i})$ once, and then multiply $({\bf B_2} \cdot \nabla \phi_{h, j}, \phi_{h,i})$ by $B_3(t)$ at each time instance to obtain  $({\bf B} \cdot \nabla \phi_{h, j}, \phi_{h,i})$. For a fixed i, there are only finite j such that $(\nabla \phi_{h, j}, \phi_{h,i})$ is not 0. This means the computational complexity of the multiplication by $B_3(t)$ is $\mathcal{O}(N_g)$. Therefore, the computational complexity for building $({\bf B} \cdot \nabla \phi_{h, j}, \phi_{h,i})$ is $\mathcal{O}(N_g)$ at each time instance. Similarly, we have that the computational complexity for  $(c\phi_{h, j}, \phi_{h,i})$ and $(f, \phi_{h,i})$ is also $\mathcal{O}(N_g)$. Therefore, the computational complexity for building the linear system (\ref{FEMform}) is $\mathcal{O}(N_g)$. We then see the cost for solving the discretized linear system. For solution of (\ref{FEMform}), since ${\bf A}_h^k$ is sparse, there are many solvers \cite{Saad1986GMRESAG, petsc} which can deal with it at a cost of $\mathcal{O}(N_g)$. Hence, the computational cost at each time instance is $\mathcal{O}(N_g)$, and the total cost during $t\in[0, T]$ is $\mathcal{O}(N_g)\times \frac{T}{\delta t}$.

We then turn to see the computational complexity for the standard POD method. When $t\in[0, T_0]$, we need to discrete the system (\ref{semi}) in the finite element space and then solve the discretized system (\ref{FEMform}). By the analysis for the standard finite element discretization, we know that the cost is $\mathcal{O}(N_g)$ at each time instance. For the construction of the POD modes, we need to perform SVD on ${\bf U}_{h}\in \mathbb{R}^{N_g\times n_s}$, where $n_{s}=\lfloor\frac{T_{0}}{\delta t\cdot\delta M}\rfloor$ as we mentioned before. Since only the left singular vectors and the singular values need to be calculated, the computational complexity is $\mathcal{O}(n_s^2 N_g)$\cite{golub13}, where $n_s^2\ll N_g$. When $t > T_0$, we need to project (\ref{semi}) on the POD subspace and then solve the discretized system (\ref{PODform}). The cost for building the linear system is $\mathcal{O}(N_g)$ at each time instance. Since ${\bf B}$, $c$ and $f$ are separable in time and space, as we mentioned before, we only need to calculate $({\bf B_1} \cdot \nabla \psi_{h, j}, \psi_{h,i})$ and $({\bf B_2} \cdot \nabla \psi_{h, j}, \psi_{h,i})$ once, which costs $\mathcal{O}(N_g)$, and then multiply it by $B_3(t)$ at each time instance to obtain  $({\bf B} \cdot \nabla \psi_{h, j}, \psi_{h,i})$, where $i, j=1,2,\cdots,m$. Since matrices composed by $({\bf B_2} \cdot \nabla \psi_{h, j}, \psi_{h,i})$ and $(c\psi_{h, j}, \psi_{h,i})$ are dense, the complexity of building the linear system (\ref{PODform}) is $\mathcal{O}(m^2)$ at each time instance. For solving the discretize system (\ref{PODform}), since ${\bf \bar A}_h^k$ is a small dense matrix, we usually use some direct method to solve it, which costs $\mathcal{O}(m^3)$. Therefore the total computational cost is $\mathcal{O}(N_g)\times (\frac{T_0}{\delta t}+1) + \mathcal{O}(n_s^2 N_g) + \mathcal{O}(m^3)\times \frac{T-T_0}{\delta t}$.

We now see the computational cost for the adaptive POD method. Except for the costs same as those for the standard POD method, some additional costs are needed for steps \textbf{Estimate}, \textbf{Mark} and \textbf{Update}. For the step \textbf{ Estimate}, different methods have different ways to design the error indicator. We denote the computational cost for this part in each time instance as $t_{\rm est}$. The cost for the step \textbf{Mark} can be neglected. The main cost for the step \textbf{Update} is that for obtaining the standard finite element approximations on interval $[t, t+\delta T]$. From the analysis for the standard finite element method, we know that the cost for this part is $\mathcal{O}(N_g)$ at each time instance. Let $n_A$ denote the number of updates for POD modes. Therefore, the total computational cost for SVD is $\mathcal{O}(n_s^2 N_g)\times 2n_A$. Since the number of POD modes will increase as the updating continues, we denote $m_A$ average number of POD modes in the adaptive POD method. By the analysis for the standard POD method, the cost for building the POD linear system is $\mathcal{O}(N_g)\times n_A + \mathcal{O}(m_A^2)\times \frac{T-T_0- n_A \delta T}{\delta t}$, and the cost for solving the discretized system at each time instance is $\mathcal{O}(m_A^3)$.  Therefore, the total computational cost is  $\mathcal{O}(N_g)\times (\frac{T_0 + n_A \cdot\delta T}{\delta t}+ n_A + 1)+ \mathcal{O}(n_s^2 N_g)\times (2n_A +1) + \mathcal{O}(m_A^3)\times \frac{T-T_0- n_A \delta T}{\delta t} + t_{\rm est}\times\frac{T-T_0- n_A \delta T}{\delta t}$.

We summarize the total computational cost for each method in Table \ref{Cost}. Usually, we have $T_0 \ll T$, $m^3 \ll N_g$ and $m_A^3 \ll N_g$. Therefore, we can see that the POD type methods usually cost less CPU time than the standard finite element method. Then we focus on the cost for the adaptive POD method. If an error indicator is cheaper, the term $t_{\rm est}$ will be less. If an error indicator is more sensitive, it will require fewer updates $n_A$ to achieve the same accuracy, then it may decrease the degree of freedom $m_A$ at the same time. Therefore, the construction of the error indicator plays an important role in reducing the cost for the adaptive POD method.
\begin{center}
    \begin{table}[H]
        \caption{Total complexity of different methods}\label{Cost}\centering
        \small{
            \begin{tabular}{|c|c|}
                \hline
                Method                             & Complexity                                                                                                      \\
                \hline
                standard finite element            & $\mathcal{O}(N_g)\times \frac{T}{\delta t}$                                                                               \\
                                \hline
                standard POD method                & $\mathcal{O}(N_g)\times (\frac{T_0}{\delta t}+1) + \mathcal{O}(n_s^2 N_g) + \mathcal{O}(m^3)\times \frac{T-T_0}{\delta t}$                    \\
                \hline
                \multirow{2}*{adaptive POD method} & $\mathcal{O}(N_g)\times (\frac{T_0 + n_A \cdot\delta T}{\delta t}+ n_A + 1)+ \mathcal{O}(n_s^2 N_g)\times (2n_A +1)$                \\
                ~                                  & + $\mathcal{O}(m_A^3)\times \frac{T-T_0- n_A \delta T}{\delta t} + t_{\rm est}\times\frac{T-T_0- n_A \delta T}{\delta t}$ \\
                \hline
            \end{tabular}
        }

    \end{table}
\end{center}

\subsection{Typical existing adaptive POD methods}
There are some existing works on adaptive POD methods \cite{Dai2017, Rapun2015, Kuang2020}. The main difference between different adaptive POD methods lies in the construction of the error indicator. Here, we introduce two typical methods, one is the residual based adaptive POD method \cite{Dai2017, Rapun2015}, the other is the two-grid based adaptive POD method \cite{Kuang2020}.

For the residual based adaptive POD method, the residual is used to construct the error indicator. In detail, the error indicator $\eta_{k}$ at time instance $t= k\delta t$ is defined as
\begin{align}
    \eta_{k} = \frac{\Vert{\bf A}_{h}^{k}{\bf \widetilde{R}} {\bf{u}}^{k}_{h, \rm{POD}} - {\bf b}_{h}^{k} - {\bf C}_{h}{\bf \widetilde{R}} {\bf{u}}^{k-1}_{h, \rm{POD}}\Vert_2}{\Vert{\bf b}_{h}^{k} + {\bf C}_{h}{\bf \widetilde{R}} {\bf{u}}^{k-1}_{h, \rm{POD}}\Vert_2}.
\end{align}
We see that the computational cost $t_{\rm est}$ for $\eta_k$ is $\mathcal{O}(N_g)$, the total computational cost for all steps \textbf{Estimate} is $\mathcal{O}(N_g)\times \frac{T-T_0- n_A\cdot \delta T}{\delta t}$.

For the two-grid based adaptive POD method, two finite element spaces are constructed, the coarse finite element space $V_{H}$ and the fine finite element space $V_{h}$. The fine finite element space is used to construct the POD modes, while the coarse finite element space is used to design the error indicator. Let $\Delta t$ be the coarse time step. The error indicator $\eta_{k}$ at time instance $t= k\Delta t$ is constructed by the approximations in the coarse finite element space $V_{H}$, that is, 
\begin{align}\label{TGAind}
    \eta_{k} = \frac{\Vert{u}_{H}^{k} - {u}_{H, \rm{POD}}^{k}\Vert_2}{\Vert u_{H}^{k}\Vert_2}.
\end{align}
Here $u_{H}^{k}$, $u_{H, \rm{POD}}^{k}$ are the standard finite element approximation and the adaptive POD approximation, respectively. For a given $\eta_{0}$, if $\eta_{k} > \eta_{0}$, the time instance $t= k\Delta t$ will be picked out.

Let $N_G$ denote the degree of freedom of the coarse space. According to the complexity analysis above, the cost for obtaining the finite element approximations in coarse space is $\mathcal{O}(N_G)\times \frac{T}{\Delta t}$, and the computational cost for obtaining the adaptive POD approximations in coarse space is $\mathcal{O}(N_G)\times (\frac{T_0 + n_A \cdot\delta T}{\Delta t}+ n_A + 1)+ \mathcal{O}(n_s^2 N_G)\times(2n_A+1) + \mathcal{O}(m_A^3)\times \frac{T-T_0- n_A \delta T}{\Delta t}$. Therefore, the computational cost for \textbf{Estimate} is $\mathcal{O}(N_G)\times (\frac{T+T_0 + n_A \cdot\delta T}{\Delta t}+ n_A + 1)+ \mathcal{O}(n_s^2 N_G)\times(2n_A+1)+ \mathcal{O}(m_A^3)\times \frac{T-T_0- n_A \delta T}{\Delta t}$. Since $N_G\ll N_g$, and $\Delta t\ll\delta t$, the cost for step \textbf{Estimate} in the two-grid based adaptive POD method is usually cheaper than that for step \textbf{Estimate} in the residual based adaptive POD method.

\section{Augmented subspace based adaptive POD method}
As we refer to above, the main difference between different adaptive POD methods is the construction of the error indicator. In this section, we introduce a new approach for developing some new adaptive POD methods. 

\subsection{General framework of the augmented subspace based adaptive POD method}
 The main idea of our new approach is to add some auxiliary modes to the current POD subspace to construct an augmented subspace, and then use the gap between the approximation obtained in the augmented subspace and that obtained in the original POD subspace to develop an error indicator.

Recall that the POD subspace is denoted as $V_{h, \rm{POD}} = \text{span}\{\psi_{h,1}, \cdots,\psi_{h,m}\}$. Suppose $\psi_{h, m+1},\cdots,$ $\psi_{h, m+r}$ be some modes which are normalized and orthogonal against each other, and orthogonal against $V_{h, \rm{POD}}$. Then, we augment the subspace $V_{h, \rm{POD}}$ by $\widetilde{V}_{h, \rm{POD}} = V_{h, \rm{POD}}\oplus \text{span}\{\psi_{h,m+1}, \cdots , \psi_{h, m+r}\}$. Next we design the error indicator $\eta_k$ at time instance $t=k\delta t$.

At time instance $t =( k-1)\delta t$, the POD approximation in the subspace $V_{h, \rm{POD}}$ can be expressed as
\begin{align}\label{uk-1}
    u^{k-1}_{h, \rm{POD}}(x, y, z)=\mathop{\sum}^{m}\limits_{i=1} \alpha_{h,i}^{k-1} \psi_{h,i}(x, y, z).
\end{align}
The approximation in the augmented subspace $\widetilde{V}_{h, \rm{POD}}$ at time instance $t = k\delta t$ can be expressed as
\begin{align}\label{uk}
    \widetilde u^{k}_{h, \rm{POD}}(x, y, z)=\mathop{\sum}^{m+r}\limits_{i=1} \widetilde{\alpha}_{h,i}^{k} \psi_{h,i}(x, y, z).
\end{align}
Inserting  (\ref{uk-1}) and (\ref{uk}) into (\ref{semi}), and setting $v=\psi_{h,j}, j=1, 2, \cdots, m+r$, respectively, we get
\begin{align}\label{detail1}
    \left( \sum^{m+r}_{i=1} \widetilde{\alpha}_{h,i}^{k} \psi_{h,i} -\sum^{m}_{i=1} \alpha_{h,i}^{k-1} \psi_{h,i}, \psi_{h, j}\right) + \delta t a(t_{k}; \sum^{m+r}_{i=1} \widetilde{\alpha}_{h,i}^{k} \psi_{h,i}, \psi_{h, j}) = \delta t(f_{k}, \psi_{h, j}).
\end{align}
The equation (\ref{detail1}) can be rewritten as
\begin{align}\label{detail2}
    ( \sum^{m+r}_{i=1} \widetilde{\alpha}_{h,i}^{k} \psi_{h,i}, \psi_{h, j})+\delta t a(t_{k}; \sum^{m+r}_{i=1} \widetilde{\alpha}_{h,i}^{k} \psi_{h,i}, \psi_{h, j}) = \delta t(f_{k}, \psi_{h, j}) + (\sum^{m}_{i=1}\alpha_{h,i}^{k-1} \psi_{h,i}, \psi_{h,j}).
\end{align}
Define
\begin{align*}
     & {\bf \widetilde A}^{k}_{h, ij}=(\psi_{h, j }, \psi_{h,i})+\delta t a(t_{k}; \psi_{h, j}, \psi_{h,i}),                                                                                              \\
     & {\bf \widetilde b}^{k}_{h}=\delta t ((f_{k}, \psi_{h,1}),\cdots,(f_{k},\psi_{h,m+r}))^{T}, {\bf \widetilde C}_{h, ij}=(\psi_{h, j}, \psi_{h,i}),                                                   \\
     & {\bf u}^{k-1}_{h, \rm{POD}}=(\alpha_{h,1}^{k-1}, \cdots, \alpha_{h,m}^{k-1})^{T}, {\bf \widetilde u}^{k}_{h, \rm{POD}}=(\widetilde{\alpha}_{h,1}^{k}, \cdots, \widetilde{\alpha}_{h,m+r}^{k})^{T}.
\end{align*}
Then we obtain the following algebraic system from (\ref{detail2})
\begin{align}\label{Algebra}
    {\bf \widetilde A}^{k}_{h} {\bf \widetilde u}^{k}_{h, \rm{POD}} = {\bf \widetilde b}^{k}_{h} +{\bf \widetilde C}_{h} {\bf u}^{k-1}_{h, \rm{POD}}.
\end{align}
We define the error indicator $\eta_{k}$ at time instance $t = k\delta t$ as
\begin{align}\label{Etak}
    \eta_{k}=\dfrac{\Vert {\widetilde u}^{k}_{h, \rm{POD}} -  u^{k}_{h, \rm{POD}}\Vert_{2}}{\Vert{\widetilde u}^{k}_{h, \rm{POD}}\Vert _2}.
\end{align}

By some formally analysis, we can see it is reasonable to use $\eta_k$ defined in (\ref{Etak}) as an error indicator. It is obvious that
\begin{align}\label{sum}
    u^{k}_{h, \rm{POD}}-\widetilde u_{h, \rm{POD}}^{k}= u^{k}_{h, \rm{POD}}-u_{h}^{k} + u_{h}^{k}- \widetilde u_{h, \rm{POD}}^{k}.
\end{align}
Since $V_{h, \rm{POD}}\subset \widetilde{V}_{h, \rm{POD}} \subset V_h$, we have
\begin{align}\label{leq}
    \Vert u_{h}^{k}-\widetilde u_{h, \rm{POD}}^{k}\Vert_{2}\leq\Vert u_{h}^{k}-u_{h, \rm{POD}}^{k}\Vert_{2}.
\end{align}
From (\ref{sum}) and (\ref{leq}), we easily obtain
\begin{align}
    \frac{1}{2}\Vert u^k_{h, \rm{POD}}-\widetilde u_{h, \rm{POD}}^{k}\Vert_{2}\leq \Vert u_{h}^{k}- u_{h, \rm{POD}}^{k}\Vert_{2}.
\end{align}
In further, if there exists $0<\zeta<1$, s.t. $\Vert u_{h}^{k}- \widetilde{u}_{h, \rm{POD}}^{k}\Vert_{2}\leq \zeta  \Vert u^k_{h}- u_{h, \rm{POD}}^{k}\Vert_{2}$, then from (\ref{sum}) we have

\begin{align*}
    \| u_h^k - u_{h, \rm{POD}}^k\|_2 \leq \| u_h^k -\widetilde{u}_{h, \rm{POD}}^k\|_2 +
    \|   u_{h, \rm{POD}}^k - \widetilde{u}_{h, \rm{POD}}^k\|_2
    \leq \zeta \| u_h^k - u_{h,\rm{POD}}^k\|_2 + \|   u_{h, \rm{POD}}^k - \widetilde{u}_{h, \rm{POD}}^k\|_2 ,
\end{align*}
from which we have
\begin{align}
    \| u_h^k - u_{h, \rm{POD}}^k\|_2 \leq  \frac{1}{1-\zeta} \|   u_{h, \rm{POD}}^k - \widetilde{u}_{h, \rm{POD}}^k\|_2 .
\end{align}
Therefore,
\begin{align}\label{inequ}
    \frac{1}{2}\Vert u^k_{h, \rm{POD}}-\widetilde u_{h, \rm{POD}}^{k}\Vert_{2}\leq \Vert u_{h}^{k}- u_{h, \rm{POD}}^{k}\Vert_{2}\leq \frac{1}{1-\zeta} \Vert u^k_{h, \rm{POD}}- \widetilde u_{h, \rm{POD}}^{k}\Vert_{2}.
\end{align}
The effect of the error indicator lies closely to $\zeta$ in (\ref{inequ}).

Apply the new error indicator defined in (\ref{Etak}) into step 8 of Algorithm \ref{APODAlg}, we then obtain the general framework of our augmented subspace based adaptive POD method, as shown in Algorithm \ref{Sub-APOD}.
\begin{algorithm}[H]
    \caption{General framework of the augmented subspace based adaptive POD method}\label{Sub-APOD}
    \begin{algorithmic}[1]
        \State Given $\delta t$, $T_{0}$, $\delta T$, $T$, $\gamma _{1}$, $\gamma_{2}$, $\gamma _{3}$, $\delta M$ and the mesh $\mathcal{T}_{h}$.
        \State Discretize (\ref{semi})  in $V_{h}$ on interval $[0, T_{0}]$ and obtain ${\bf u}^{k}_{h}$, $\forall k \in [0, \lfloor T_{0}/\delta t\rfloor]$, then obtain snapshots ${\bf U}_{h}$ at different times $t_{0}$, $t_{\delta M}$, $\cdots$, $t_{n_{s}\cdot \delta M}$.
        \State Construct POD modes $\Psi_{h}=\{\psi_{h,1}, \cdots,\psi_{h,m}\}$  by POD\_Mode(${\bf U}_{h}, \gamma_1, \Phi_{h}, m, \Psi_{h}$).
        \State $t = T_{0}$, $k = \frac{T_{0}}{\delta t}$.
        \While {$t \leq T$}
        \State  $t = t+\delta t$, $k = k+1$.
        \State  Discretize (\ref{semi}) in the POD subspace $V_{h, \rm{POD}}= \text{span}\{\psi_{h,1}, \cdots,\psi_{h,m}\}$, and obtain the POD approximation ${u}^{k}_{h, \rm{POD}}$.
        \State Provide some auxiliary modes and orthogonalize them against $\Psi_{h}$, and make them normalized and orthogonal against each other, then obtain the modes $\psi_{h,m+1},\cdots,\psi_{h,m+r}$.
        \State  Discretize (\ref{semi}) in the augmented subspace $V_{h, \rm{POD}}\oplus \text{span}\{\psi_{h,m+1}, \cdots , \psi_{h, m+r}\}$, and obtain ${\widetilde u}^{k}_{h, \rm{POD}}$.
        \State Compute the error indicator $\eta_{k}$ by (\ref{Etak}).
        \If {$\eta_{k} > \eta_{0}$}
        \State  $t=t-\delta t$, $k=k-1$.
        \State Discretize (\ref{semi}) in $V_{h}$ on interval $[t, t+\delta T]$ and get snapshots ${\bf W}_{h, 1}$, then update POD modes $\Psi_{h}$  by Update\_POD\_Mode(${\bf W}_{h, 1}, \gamma_{2}, \gamma_{3}, \Phi_{h}, m, \Psi_{h}$).
        \State $k=k+ \frac{\delta T}{\delta t}$.
        \EndIf
        \State {\bf end if}
        \EndWhile
        \State{\bf end while}
    \end{algorithmic}
\end{algorithm}

\subsection{Specific augmented subspace based adaptive POD methods}
The key for constructing the auxiliary modes includes two points: one is that the auxiliary modes can not be orthogonal to the exact solution, and it is better that the angle between these auxiliary modes and the exact solution is far away from $\pi/2$, the other is that they should be cheap to be constructed. We now provide two specific methods for obtaining the auxiliary modes. For simplicity, we only consider the case of $r=1$.
\subsubsection{Random vector based augmented subspace}
One of the simplest methods for augmenting the current POD subspace is to use a randomly generated vector as the auxiliary mode. We denote this auxiliary mode by {\bf d}. In fact, it is basically impossible that the random vector is linearly correlated with the current POD modes. We orthogonalize {\bf d} against ${\bf \widetilde{R}}$ and normalize it, and still denote it as {\bf d} in our following analysis.

Using the similar procedure for deducing the linear system (\ref{Algebra}), we obtain
\begin{align}\label{Randform}
    \begin{bmatrix}
        {\bf \widetilde{R}}^{T}{\bf A}^{k}_{h}{\bf \widetilde{R}} & {\bf \widetilde{R}}^{T}{\bf A}^{k}_{h}{\bf d} \\
        {\bf d}^{T}{\bf A}^{k}_{h}{\bf \widetilde{R}}             & {\bf d}^{T}{\bf A}^{k}_{h}{\bf d}
    \end{bmatrix}
    {\bf \widetilde{u}}^{k}_{h, \rm{POD}}
    =\begin{bmatrix}
        {\bf \widetilde{R}}^{T}{\bf b}^{k}_{h} \\
        {\bf d}^{T}{\bf b}^{k}_{h}
    \end{bmatrix}
    +
    \begin{bmatrix}
        {\bf \widetilde{R}}^{T}{\bf C}_{h}{\bf \widetilde{R}} \\
        {\bf d}^{T}{\bf C}_{h}{\bf \widetilde{R}}
    \end{bmatrix}
    {\bf u}^{k-1}_{h, \rm{POD}}.
\end{align}
The computational cost of the step \textbf{Estimate} is as follows. The implementation of orthogonalization requires vector-vector products $m_A$ times, where $m_A$ is the number of vectors in ${\bf \widetilde{R}}$, and normalization  requires vector-vector products only once. Therefore, the computational cost of orthogonal normalization is  $\mathcal{O}(N_g)$. Since the matrix in (\ref{Randform}) can be assembled by the previously calculated values for the POD linear system in (\ref{PODform}) except for the terms which contain ${\bf d}$, we only require to compute the elements of the last row and the last column of the matrix. Due to ${\bf A}^{k}_{h}$ and ${\bf C}_{h}$ are sparse, the cost for building the augmented linear system is $\mathcal{O}(N_g)$. Since the size of the linear system in (\ref{Randform}) is only 1 more than that of POD linear system, similar to the analysis for the adaptive POD method, the cost for computing ${\bf \widetilde{u}}^{k}_{h, \rm{POD}}$ is $\mathcal{O}(m_A^3)$. The computational cost for $\eta_k$ in (\ref{Etak}) is $\mathcal{O}(N_g)$. Therefore, the cost for step \textbf{Estimate} at each time instance is $\mathcal{O}(N_g)+\mathcal{O}(m_A^3)$.

We then see the angle $\theta$ between this auxiliary mode and the exact solution. Here, we take the equation (\ref{equ}) with the advection being the ABC flow for case of $\epsilon = 0.01$ as an example. We use the standard finite element approximation as the exact solution. The the cosine of $\theta$ $\cos\theta$ obtained at each time instance $t\in[0,100]$ are $10^{-5}\sim10^{-3}$.

We can see that the $\cos\theta$ is very close to zero, which means the augmented mode is almost orthogonal to the exact solution. Therefore, although using a random vector as an auxiliary mode is cheap, the effect is not as good as we expect. It is not a recommended method to augment the subspace. 
\subsubsection{Coarse-grid approximations based augmented subspace}
We see from \cite{Kuang2020} that the  solution obtained in the coarse finite element space is a good approximation for the solution obtained in the fine finite element space. Moreover, the computational cost for obtaining the solution approximation in the coarse finite element space is far less than that for obtaining the solution approximation in the fine finite element space. Therefore, here, we consider to use the approximated solution obtained in the coarse finite element space as the auxiliary mode to construct the augmented subspace.

We denote $u^l_{H}$ the finite element approximation in coarse finite element space $V_H$ at each time instance $t=l\Delta t$, where $ l=0, 1,\cdots, {\lfloor\frac{T}{\Delta t}\rfloor}$ and  $\Delta t$ also denotes the coarse time step. We denote the interpolation of $u_H^{l}$ in the fine finite element space by $u_{H, I}^{l}$. Then, we orthogonalize ${\bf u}^l_{H,I}$ against ${\bf \widetilde{R}}$ and normalize it, and denote it as $\psi^l_{h, m+1}$. We augment the subspace at $t=l\Delta t$ by
\begin{align}
    {\widetilde V}^l_{h, \rm{POD}}:=V_{h, \rm{POD}}\oplus \text{span}\{\psi^l_{h,m+1}\}.
\end{align}
Since $\delta t \ll\Delta t$, we set $k$ by $k = \frac{\Delta t}{\delta t} l$. The approximation in the augmented subspace ${\widetilde V}^l_{h, \rm{POD}}$ at $t=l\Delta t = k\delta t$ can be expressed as
\begin{align}\label{tildeuk}
    \widetilde u^{k}_{h, \rm{POD}}(x, y, z)=\mathop{\sum}^{m}\limits_{i=1} \widetilde{\alpha}_{h,i}^{k} \psi_{h,i}(x, y, z) + \widetilde{\alpha}_{h,m+1}^{k} \psi^{l}_{h,m+1}(x, y, z).
\end{align}
We know that the POD approximation in the subspace $V_{h, \rm{POD}}$ at $t =( k-1)\delta t$ can be expressed as
\begin{align}\label{tildeuk-1}
    u^{k-1}_{h, \rm{POD}}(x, y, z)=\mathop{\sum}^{m}\limits_{i=1} \alpha_{h,i}^{k-1} \psi_{h,i}(x, y, z).
\end{align}
Similar to the derivation in Section 3.1, we insert (\ref{tildeuk-1}) and (\ref{tildeuk}) into (\ref{semi}), then set $v=\psi_{h,j}, j=1, 2, \cdots, m$ and $v=\psi^{l}_{h,m+1}$, respectively. We obtain the following the equation
\begin{equation}\label{cosfull}
\begin{split}
    ( \sum^{m}_{i=1} \widetilde{\alpha}_{h,i}^{k} \psi_{h,i} + \widetilde{\alpha}_{h,m+1}^{k} \psi^{l}_{h,m+1}, \psi_{h, j}) & +\delta t a(t_{k}; \sum^{m}_{i=1} \widetilde{\alpha}_{h,i}^{k} \psi_{h,i} + \widetilde{\alpha}_{h,m+1}^{k} \psi^{l}_{h,m+1}, \psi_{h, j}) \\
                                                                                                                             & = \delta t(f_{k}, \psi_{h, j}) + (\sum^{m}_{i=1}\alpha_{h,i}^{k-1} \psi_{h,i}, \psi_{h, j}).
\end{split}
\end{equation}

In the practical operation, we only need to build an augmented subspace ${\widetilde V}^l_{h, \rm{POD}}$ to calculate the error indicator $\eta_l$ each coarse time step $\Delta t$. We compute the error indicator $\eta_l$ by
 \begin{align}\label{CosEtak}
              \eta_{l}=\dfrac{\Vert {\bf \widetilde u}^{k}_{h, \rm{POD}} -  {\bf u}^{k}_{h, \rm{POD}}\Vert_{2}}{\Vert{\bf \widetilde u}^{k}_{h, \rm{POD}}\Vert _2}.
          \end{align}
For the convenience of the following discussion, we summarize the process for computing the error indicator as routine Error\_Indicator($l,{\bf A}^{k}_{h}, {\bf b}^{k}_{h}, {\bf C}_{h}, {\bf \bar A}^{k}_{h}, {\bf \bar b}^{k}_{h}, {\bf \bar C}_{h}, {\bf u}^{k-1}_{h, \rm{POD}}, {\bf u}^l_{H,I}, \Psi_{h}, \eta_l$) in Algorithm \ref{ErrInd-Alg}.
\begin{algorithm}[H]
    \caption{Error\_Indicator($l,{\bf A}^{k}_{h}, {\bf b}^{k}_{h}, {\bf C}_{h}, {\bf \bar A}^{k}_{h}, {\bf \bar b}^{k}_{h}, {\bf \bar C}_{h}, {\bf u}^{k-1}_{h, \rm{POD}}, {\bf u}^l_{H,I}, \Psi_{h}, \eta_l$)}\label{ErrInd-Alg}
    {\bf Input:} $l,{\bf A}^{k}_{h}, {\bf b}^{k}_{h}, {\bf C}_{h}, {\bf \bar A}^{k}_{h}, {\bf \bar b}^{k}_{h}, {\bf \bar C}_{h}, {\bf u}^{k-1}_{h, \rm{POD}}$, ${\bf u}^l_{H,I}$ and $\Psi_{h}$, $\Psi_{h} = \Phi_{h}{\bf \widetilde R}$.\\
    {\bf Output:} $\eta_l$.
    \begin{algorithmic}[1]
        \State  Orthogonalize ${\bf u}^l_{H,I}$ against ${\bf \widetilde{R}}$ and normalize it, then obtain the auxiliary vector ${\bf d}_l$ at $t=l\Delta t$.
        \State Compute ${\bf \widetilde{u}}^{k}_{h, \rm{POD}}$ in (\ref{cosfull}) at $t=k\delta t= l\Delta t$ by
          \begin{align}\label{Cosform}
              \begin{bmatrix}
                  {\bf \bar A}^{k}_{h}  & {\bf \widetilde{R}}^{T}{\bf A}^{k}_{h}{\bf d}_l \\
                  {\bf d}_l^{T}{\bf A}^{k}_{h}{\bf \widetilde{R}}           & {\bf d}_l^{T}{\bf A}^{k}_{h}{\bf d}_l
              \end{bmatrix}
              {\bf \widetilde{u}}^{k}_{h, \rm{POD}}
              =\begin{bmatrix}
                 {\bf \bar b}^{k}_{h} \\
                  {\bf d}_l^{T}{\bf b}^{k}_{h}
              \end{bmatrix}
              +
              \begin{bmatrix}
                   {\bf \bar C}_{h}\\
                  {\bf d}_l^{T}{\bf C}_{h}{\bf \widetilde{R}}
              \end{bmatrix}
              {\bf u}^{k-1}_{h, \rm{POD}}.
          \end{align}
        \State Obtain the error indicator $\eta_l$ at time instance $t= l\Delta t$ by (\ref{CosEtak}).
       
    \end{algorithmic}
\end{algorithm}
Apply Algorithm \ref{ErrInd-Alg} into step 10 of Algorithm \ref{Sub-APOD}, we then obtain the augmented subspace based adaptive POD method with coarse-grid approximations, as shown in Algorithm \ref{Aug-APOD}.
\begin{algorithm}[H]
    \caption{Augmented subspace based adaptive POD method with coarse-grid approximations}\label{Aug-APOD}
    \begin{algorithmic}[1]
        \State Given $\delta t$, $\Delta t$, $T_{0}$, $\delta T$, $T$, $\gamma _{1}$, $\gamma_{2}$, $\gamma _{3}$, $\delta M$ and the mesh $\mathcal{T}_{h}$, $\mathcal{T}_{H}$.
        \State Discretize (\ref{semi}) in $V_{H}$ on interval $[0, T]$, and obtain the approximations $\{u_H^l\}$, $\forall l\in[0, {\lfloor\frac{T}{\Delta t}\rfloor}]$.
        \State Interpolate $\{u_H^l\}$ to the fine finite element space, then obtain the interpolations $\{u_{H,I}^l\}$.
        \State Discretize (\ref{semi})  in $V_{h}$ on interval $[0, T_{0}]$ and obtain ${\bf A}^{k}_{h}, {\bf b}^{k}_{h}, {\bf C}_{h}, {\bf u}^{k}_{h}$, $\forall k \in [0, {\lfloor T_{0}/\delta t\rfloor}]$, then obtain snapshots ${\bf U}_{h}$ at different times $t_{0}$, $t_{\delta M}$, $\cdots$, $t_{n_{s}\cdot \delta M}$.
        \State Construct POD modes $\Psi_{h}=\{\psi_{h, 1}, \cdots, \psi_{h, m}\}$  by POD\_Mode(${\bf U}_{h}, \gamma_1, \Phi_{h}, m, \Psi_{h}$).
        \State $t = T_{0}$, $k=\frac{T_{0}}{\delta t}$, $w=\frac{\Delta t}{\delta t}$.
        \While {$t \leq T$}
        \State  $t = t+\delta t$, $k = k+1$.
        \State  Discretize (\ref{semi}) in the subspace $ V_{h, \rm{POD}}=\text{span}\{\psi_{h, 1}, \cdots, \psi_{h, m}\}$, then obtain ${\bf \bar A}^{k}_{h}, {\bf \bar b}^{k}_{h}, {\bf \bar C}_{h}$ and ${\bf u}^{k}_{h, \rm{POD}}$.
        \If {$k\% w=0$}
        \State $l=\frac{k}{w}$.
        \State Compute error indicator $\eta_l$ by Error\_Indicator($l,{\bf A}^{k}_{h}, {\bf b}^{k}_{h}, {\bf C}_{h}, {\bf \bar A}^{k}_{h}, {\bf \bar b}^{k}_{h}, {\bf \bar C}_{h}, {\bf u}^{k-1}_{h, \rm{POD}}, {\bf u}^l_{H,I}, \Psi_{h}, \eta_l$).
        \If {$\eta_{l} > \eta_{0}$}
        \State  $t=t-\delta t$, $k=k-1$.
        \State Discretize (\ref{semi})  in $V_{h}$ on interval $[t, t+\delta T]$ to get ${\bf u}^{k+i}_{h}, i=1, \cdots,\frac{\delta T}{\delta t}$, then obtain snapshots ${\bf W}_{h, 1}$.
        \State Update POD modes $\Psi_{h}$ by Update\_POD\_Mode(${\bf W}_{h, 1}, \gamma_{2}, \gamma_{3}, \Phi_{h}, m, \Psi_{h}$). $k=k+ \frac{\delta T}{\delta t}$.
        \EndIf
        \State {\bf end if}
        \EndIf
        \State {\bf end if}
        \EndWhile
        \State{\bf end while}
    \end{algorithmic}
\end{algorithm}

Next we analyze the computational complexity of the step \textbf{Estimate} for this new method. As we mentioned before, the cost for obtaining the finite element approximation at each time instance in coarse finite element space is $\mathcal{O}(N_G)$. The cost for interpolation of a function in the coarse finite element space in the fine finite element space is $\mathcal{O}(N_g)$. Similar to the analysis in Section 3.2.1, the computational cost of orthogonal normalization is  $\mathcal{O}(N_g)$.
 Since the matrix in (\ref{Cosform}) can be assembled by the previously calculated values for the POD linear system in (\ref{PODform}) except for the terms which contain ${\bf d}_{l}$, we only require to compute the elements of the last row and the last column of the matrix. Due to ${\bf A}^{k}_{h}$ and ${\bf C}_{h}$ are sparse, the cost for building the augmented linear system is $\mathcal{O}(N_g)$. Since the size of the linear system in (\ref{Cosform}) is only 1 more than that of POD linear system, similar to the analysis for the adaptive POD method, the cost for computing ${\bf \widetilde{u}}^{k}_{h, \rm{POD}}$ is $\mathcal{O}(m_A^3)$. It is obvious that the computational cost for $\eta_l$ in (\ref{CosEtak}) is $\mathcal{O}(N_g)$.  Therefore, the total cost for \textbf{Estimate} is $(\mathcal{O}(N_g)+\mathcal{O}(m_A^3))\times \frac{T-T_0- n_A\cdot \delta T}{\Delta t} + (\mathcal{O}(N_G)+ \mathcal{O}(N_g))\times \frac{T}{\Delta t}$. As we know, $\Delta t\ll \delta t$, the number of computation times for the error indicator is relatively small. Noting that $N_G\ll N_g$, we estimate the cost for the error indicator is relatively cheap. 

\section{Numerical examples}
In this section, we will use two typical fluid advection fields with chaotic streamlines, the Kolmogorov flow and the ABC flow, to show the accuracy and efficiency of our augmented subspace based adaptive POD method.

We use the standard finite element method approximation as the reference solution, and compare our new method with the standard POD method and the two-grid based adaptive POD method, respectively. Here we compute the error indicator $\eta_l$ of our new approach by (\ref{CosEtak}). The relative error of approximation at each time instance is calculated by
\begin{align}\label{error}
    \text{Error}=\dfrac{\Vert u_{h}^{k} - u_{h, *}^{k}\Vert_{2}}{\Vert u_{h}^{k}\Vert_{2}},
\end{align}
where $u_{h}^{k}$ and $u_{h,*}^{k}$ represent the finite element approximations and different types of the POD approximations at different times $t=t_{k}$, respectively. 
The numerical experiments are carried out on the high performance computers of the State Key Laboratory of Scientific and Engineering Computing, Chinese Academy of Sciences, and our code is based on the toolbox PHG \cite{phg}.

In the following discussions, we will denote the standard finite element method, the standard POD method, the two-grid based adaptive POD method and the augmented subspace based adaptive POD method with coarse-grid approximations as ``FEM'', ``POD'', ``TG-APOD'' and ``Aug-APOD'', respectively.
\subsection{Kolmogorov flow}
We consider the following advection-diffusion equation with the advection being the Kolmogorov flow \cite{Obukhov1983, Borue1996},
\begin{equation}\label{Kol_eq}
    \left\{
    \begin{aligned}
         & u_{t} - \epsilon\Delta u + {\bf B}(x, y, z, t)\cdot \nabla u = f(x, y, z, t),\quad (x, y, z)\in\Omega, t \in [0, T], \\
         & u(x, y, z, 0) = 0,                                                                                                    \\
         & u(x + 2\pi, y, z, t) = u(x, y + 2\pi, z, t) = u(x, y, z + 2\pi, t) = u(x, y, z, t),
    \end{aligned}
    \right.
\end{equation}
where
\begin{align*}
     & {\bf B} (x, y, z, t)  = (\cos(y), \cos(z), \cos(x)) + (\sin(z), \sin(x), \sin(y))\cos(t), \\
     & f(x, y, z, t)         = - \cos(y) - \sin(z) * \cos(t),                                    \\
     & \Omega = [0, 2\pi]^{3}, T = 100.
\end{align*}
We will test 6 different cases with $\epsilon=1, 0.5, 0.1, 0.05, 0.01, 0.005$, respectively.

 We first divide the $\Omega$ into 6 tetrahedrons to act as the initial mesh. Then we refine the initial mesh 23 times uniformly using bisection to obtain the fine mesh for cases of $\epsilon = 1, 0.5, 0.1, 0.05$, and refine the initial mesh 24 times to obtain the fine mesh for cases of $\epsilon = 0.01, 0.005$. We use the piecewise linear function as the finite element basis and set $\delta t = 5\times 10^{-3}$. For the POD method, we set $T_{0} = 5$, $\delta M = 20$. For the adaptive POD methods, we set $\delta T = 4$. In all the numerical experiments, we choose the parameters $\gamma_{i}(i = 1, 2, 3)$ as $\gamma_{1} =\gamma_{2} = 0.999$, $\gamma_{3} = 1.0-1.0\times10^{-8}$. For the methods TG-APOD and Aug-APOD, we refine the initial mesh 14 times to obtain the coarse mesh for cases of $\epsilon = 1, 0.5, 0.1, 0.05$, and refine the initial mesh 15 times to obtain the coarse mesh for cases of $\epsilon = 0.01, 0.005$. The time steps corresponding to the coarse finite element spaces are 0.2 and 0.125, respectively. We use 72 processors to do the simulation for cases of the fine mesh being obtained by refining the initial mesh 23 times, and 180 processors for cases of the fine mesh being obtained by refining the initial mesh 24 times. 

\begin{figure}
    \centering
    \subfloat[$\epsilon=1.$]{
        \includegraphics[width=5.5cm]{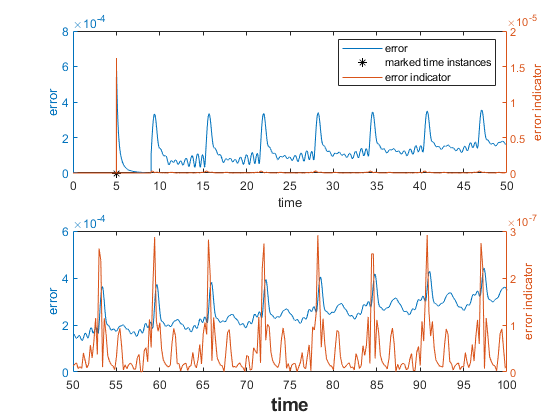}
    }
    \quad
    \subfloat[$\epsilon=0.5.$]{
        \includegraphics[width=5.5cm]{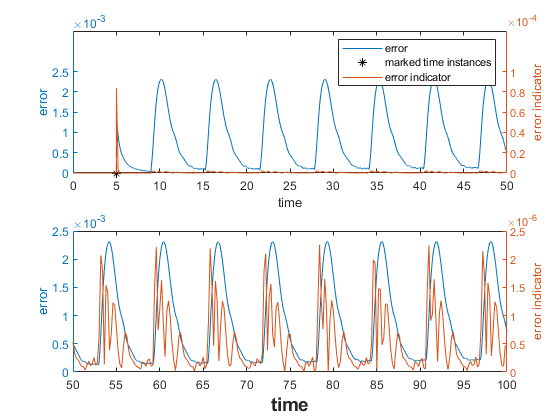}
    }
    \quad
    \subfloat[$\epsilon=0.1.$]{
        \includegraphics[width=5.5cm]{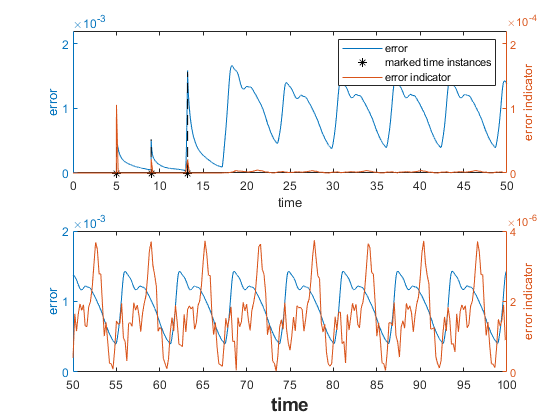}
    }
    \quad
    \subfloat[$\epsilon=0.05.$]{
        \includegraphics[width=5.5cm]{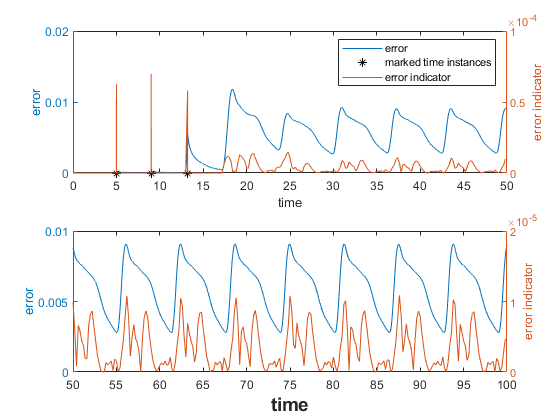}
    }
    \quad
    \subfloat[$\epsilon=0.01.$]{
        \includegraphics[width=5.5cm]{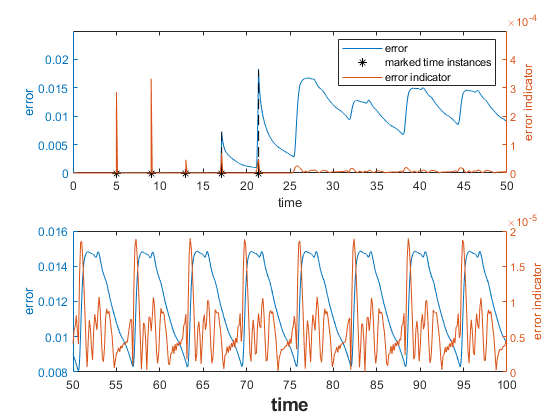}
    }
    \quad
    \subfloat[$\epsilon=0.005.$]{
        \includegraphics[width=5.5cm]{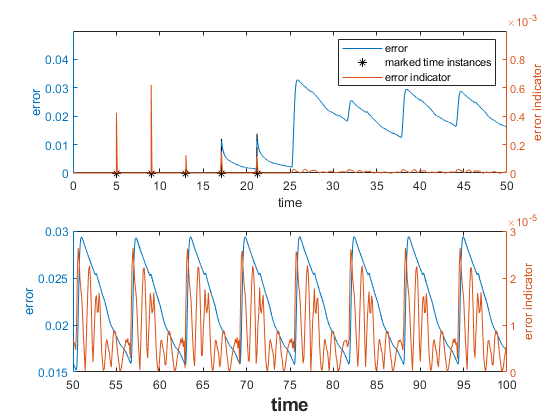}
    }
    \caption{The evolution curves of the error indicator and the error obtined by the method Aug-APOD for solution of (\ref{Kol_eq}) with $\epsilon=1,0.5,0.1,0.05,0.01,0.005$, respectively}
    \label{figKolmo}
\end{figure}

We first show the error and the error indicator obtained by our method Aug-APOD in Figure \ref{figKolmo}. The sub-figures at the top describe the variation curves of error and error indicator on time interval $[0,50]$, and the sub-figures at the bottom describe those on time interval $[50,100]$. The x-axis of each sub-figure is the time, the y-axis in the left of each sub-figure is the error which is defined by (\ref{error}), while the y-axis in the right of each sub-figure is the error indicator defined by (\ref{CosEtak}). In each sub-figure, the blue curve describes the variation of the error, the orange curve describes the variation of the error indicator, and the black star denotes the marked time instance when the POD modes are to be updated.

From each sub-figure at the bottom of Figure \ref{figKolmo}, we see that the variation of the error indicator is proportional to that of the error. Especially, the time instances when the error indicator achieves its local maximizer coincide very well with those when the error achieves its local maximizer, which means that the error indicator is very effective.
\graphicspath{{Img/Err&ErrInd/}}

We then compare our method with the methods FEM, POD and TG-APOD in Table \ref{tableKolmo}, where different methods are shown. The parameters set for the method TG-APOD and the method Aug-APOD are the same expect the threshold value $\eta_0$.

In Table \ref{tableKolmo}, `Update Times' means the number of updates for the POD modes in the adaptive POD methods,  `DOFs' means the degree of freedom, `Error' denotes the relative error for numerical solution at t=T, `Average Error' denotes the average of relative error for numerical approximation at each time instance on time interval $[0,T]$ and `Time' means the wall time for the simulation.

From Table \ref{tableKolmo}, we can first see that although the number of updates and the degrees of freedom for the POD type methods increase with the decrease of  $\epsilon$, the degrees of freedom for the POD type methods for all cases of different $\epsilon$ are much smaller than those for the standard finite element method. We also see that as the decrease of $\epsilon$, the error obtained by the standard POD method increases dramatically, which makes the results unreliable when $\epsilon$ is close to zero. This show us that the smaller the $\epsilon$, the more difficult the model to be simulated. Fortunately, results obtained by the two adaptive POD methods are still of high accuracy even for case of $\epsilon=0.005$, and the CPU time cost by the adaptive POD methods is less than one-half of that used by the standard finite element method. This shows that the adaptive POD methods behave much better than the standard POD method. We now compare the two adaptive POD methods Aug-APOD and TG-APOD. We see from Table \ref{tableKolmo} that the both two methods cost similar time with the similar degrees of freedom, but the new method Aug-APOD can achieve higher accuracy, especially for cases of $\epsilon$ being very close to zero. Specifically, for cases of $\epsilon=0.01$ and $\epsilon =0.005$, our new method can achieve higher accuracy by using less POD modes and less CPU time, and its average error can even be an order of magnitude smaller than that obtained by method TG-APOD. 

\begin{figure}[H]
    \centering
    \subfloat[$\epsilon=1.$]{
        \includegraphics[width=5.5cm]{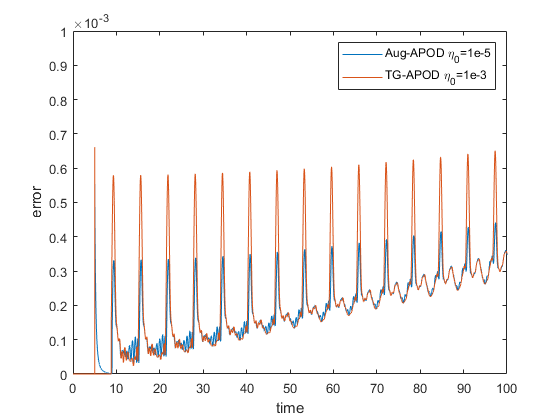}
    }
    \quad
    \subfloat[$\epsilon=0.5.$]{
        \includegraphics[width=5.5cm]{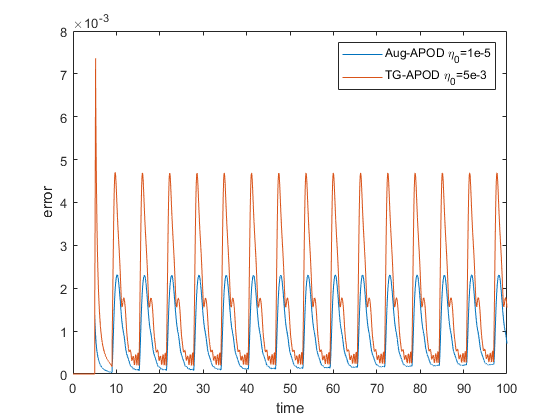}
    }
    \quad
    \subfloat[$\epsilon=0.1.$]{
        \includegraphics[width=5.5cm]{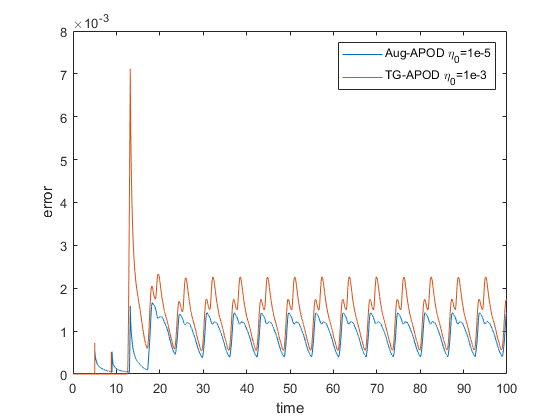}
    }
    \quad
    \subfloat[$\epsilon=0.05.$]{
        \includegraphics[width=5.5cm]{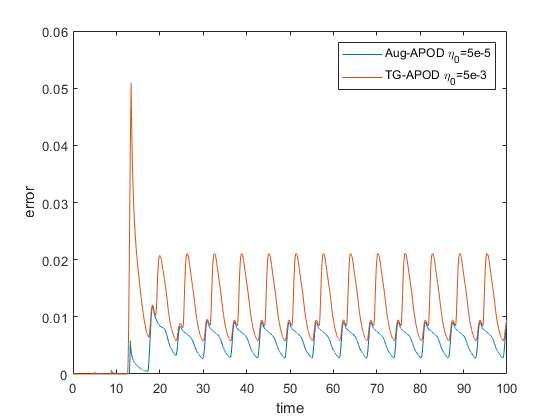}
    }
    \quad
    \subfloat[$\epsilon=0.01.$]{
        \includegraphics[width=5.5cm]{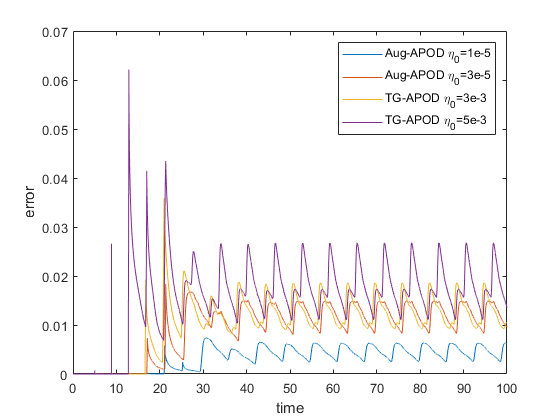}
    }
    \quad
    \subfloat[$\epsilon=0.005.$]{
        \includegraphics[width=5.5cm]{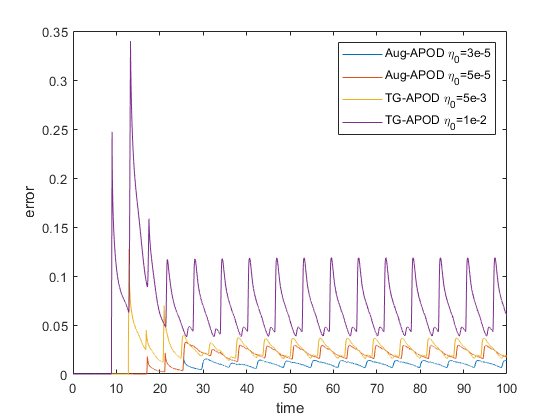}
    }
    \caption{The evolution curves of the error for solution of (\ref{Kol_eq}) with different $\epsilon$ by methods TG-APOD and Aug-APOD, respectively}
    \label{figKolmoTGCos}
\end{figure}

We have done more tests with other $\eta_0$ and different coarse meshes for cases of  $\epsilon=0.01$ and $\epsilon= 0.005$. The corresponding results are provided in ``Appendix A ''.
\begin{center}
    \begin{table}[H]
        \caption{The results of (\ref{Kol_eq}) with different $\epsilon$ obtained by methods FEM, POD, TG-APOD and Aug-APOD, respectively}\centering
        \label{tableKolmo}
        \small{
            \begin{tabular}{|cccccccc|}
                \hline
                $\epsilon$ & Method                           & $\eta_{0}$                   & Update Times & DOFs     & Error    & Average Error & Time(s)  \\
                \hline
                1          & \textcolor[rgb]{0,0,0}{FEM}      & -                            & -            & 10485760 & -        & -             & 15041.77 \\
                           & \textcolor[rgb]{0,0,0}{POD}      & -                            & -            & 19       & 0.003877 & 0.001094      & 1859.07  \\
                           & \textcolor[rgb]{0,0,0}{TG-APOD}  & \textcolor[rgb]{0,0,0}{1e-3} & 1            & 13       & 0.000350 & 0.000207      & 1829.59  \\
                           & \textcolor[rgb]{0,0,0}{Aug-APOD} & 1e-5                         & 1            & 14       & 0.000347 & 0.000184      & 1900.86  \\
                           \hline
                0.5        & \textcolor[rgb]{0,0,0}{FEM}      & -                            & -            & 10485760 & -        & -             & 14772.71 \\
                           & \textcolor[rgb]{0,0,0}{POD}      & -                            & -            & 22       & 0.002337 & 0.009928      & 1805.14  \\
                           & \textcolor[rgb]{0,0,0}{TG-APOD}  & \textcolor[rgb]{0,0,0}{5e-3} & 1            & 15       & 0.001546 & 0.001568      & 1685.59  \\
                           & \textcolor[rgb]{0,0,0}{Aug-APOD} & 1e-5                         & 1            & 16       & 0.000706 & 0.000808      & 1726.65  \\
                           \hline
                0.1        & \textcolor[rgb]{0,0,0}{FEM}      & -                            & -            & 10485760 & -        & -             & 14998.48 \\
                           & \textcolor[rgb]{0,0,0}{POD}      & -                            & -            & 29       & 0.140401 & 0.223960      & 1943.66  \\
                           & \textcolor[rgb]{0,0,0}{TG-APOD}  & \textcolor[rgb]{0,0,0}{1e-3} & 3            & 47       & 0.001719 & 0.001238      & 3154.18  \\
                           & \textcolor[rgb]{0,0,0}{Aug-APOD} & 1e-5                         & 3            & 48       & 0.001428 & 0.000839      & 3182.90  \\
\hline
                0.05       & \textcolor[rgb]{0,0,0}{FEM}      & -                            & -            & 10485760 & -        & -             & 13188.19 \\
                           & \textcolor[rgb]{0,0,0}{POD}      & -                            & -            & 34       & 0.325160 & 0.415044      & 1448.42  \\
                           & \textcolor[rgb]{0,0,0}{TG-APOD}  & \textcolor[rgb]{0,0,0}{5e-3} & 3            & 58       & 0.009429 & 0.010761      & 3350.81  \\
                           & \textcolor[rgb]{0,0,0}{Aug-APOD} & 5e-5                         & 3            & 59       & 0.008989 & 0.005068      & 3288.42  \\
                           \hline
                0.01       & \textcolor[rgb]{0,0,0}{FEM}      & -                            & -            & 16777216 & -        & -             & 11784.99 \\
                           & \textcolor[rgb]{0,0,0}{POD}      & -                            & -            & 43       & 0.919008 & 0.777626      & 1267.46  \\
                           & \textcolor[rgb]{0,0,0}{TG-APOD}  & 5e-3                         & 5            & 144      & 0.013528 & 0.015620      & 4511.11  \\
                           & \textcolor[rgb]{0,0,0}{TG-APOD}  & 3e-3                         & 6            & 165      & 0.010431 & 0.010523      & 5043.79  \\
                           & \textcolor[rgb]{0,0,0}{Aug-APOD} & 3e-5                         & 5            & 147      & 0.009588 & 0.009403      & 4505.68  \\
                           & \textcolor[rgb]{0,0,0}{Aug-APOD} & 1e-5                         & 6            & 172      & 0.006208 & 0.003351      & 5280.25  \\
                           \hline
                0.005      & \textcolor[rgb]{0,0,0}{FEM}      & -                            & -            & 16777216 & -        & -             & 12526.50 \\
                           & \textcolor[rgb]{0,0,0}{POD}      & -                            & -            & 45       & 1.010868 & 0.871304      & 1306.83  \\
                           & \textcolor[rgb]{0,0,0}{TG-APOD}  & 1e-2                         & 5            & 170      & 0.059025 & 0.037340      & 5040.39  \\
                           & \textcolor[rgb]{0,0,0}{TG-APOD}  & 5e-3                         & 6            & 200      & 0.018382 & 0.022522      & 5713.98  \\
                           & \textcolor[rgb]{0,0,0}{Aug-APOD} & 5e-5                         & 5            & 176      & 0.017100 & 0.017024      & 5084.80  \\
                           & \textcolor[rgb]{0,0,0}{Aug-APOD} & \textcolor[rgb]{0,0,0}{3e-5} & 6            & 207      & 0.013551 & 0.008223      & 5913.76  \\
                \hline
            \end{tabular}
        }
    \end{table}
\end{center}
To see more clearly, we show the variation curves of the error obtained by the two adaptive POD methods in Figure \ref{figKolmoTGCos}. The x-axis of each figure is the time, the y-axis is the relative error of numerical solution obtained by the adaptive POD methods. For the cases of $\epsilon=1, 0.5, 0.1, 0.05$, the blue curve denotes the error obtained by the method Aug-APOD, and the orange curve denotes the error obtained by the  method TG-APOD. For the cases of $\epsilon=0.01, 0.005$, the blue curve and orange curve denote the error obtained by the method Aug-APOD with different $\eta_0$, respectively, and the yellow curve and purple curve denote the error obtained by the method TG-APOD with different $\eta_0$, respectively. If we take a detail look at the error obtained by the two methods over the entire time interval, we can find that the error obtained by our method Aug-APOD is always smaller than that obtained by the method TG-APOD. Moreover, taking into account the wall time reported in Table \ref{tableKolmo}, the simulation denoted by orange curve costs the least CPU time for case of $\epsilon=0.01$. From these comparison, we can see that our new method Aug-APOD is more efficient than the method TG-APOD. 
\graphicspath{{Img/TGA&Aug_Err/}}

\subsection{Arnold-Beltrami-Childress(ABC) flow}
We then consider another advection-diffusion equations with the advection being the ABC flow, which plays an important role in fluid dynamics \cite{Dombre1986, Xin2016, Brummell2001}.
\begin{equation}\label{ABCeq}
    \left\{
    \begin{aligned}
         & u_{t} - \epsilon\Delta u + {\bf B}(x, y, z, t)\cdot \nabla u = f(x, y, z, t),\quad (x, y, z)\in\Omega, t \in [0, T], \\
         & u(x, y, z, 0) = 0,                                                                                                    \\
         & u(x + 2\pi, y, z, t) = u(x, y + 2\pi, z, t) = u(x, y, z + 2\pi, t) = u(x, y, z, t),
    \end{aligned}
    \right.
\end{equation}
where
\begin{align*}
    {\bf B} (x, y, z, t) =                 & (\sin(z + \sin wt) + \cos(y + \sin wt), \sin(x + \sin wt)    \\
                                           & + \cos(z + \sin wt), \sin(y + \sin wt) + \cos(x + \sin wt)), \\
    f(x, y, z, t)         =                & -\sin(z + \sin wt) - \cos(y + \sin wt),                      \\
    \Omega                = [0, 2\pi]^{3}, & T = 100.
\end{align*}
For this example, we also test 6 different cases with $\epsilon = 1, 0.5, 0.1, 0.05, 0.01, 0.005$, respectively. 

Similar to the example of Kolmogorov flow, we first divide the $\Omega$ into 6 tetrahedrons to act as the initial mesh. Then we refine the initial mesh 23 times uniformly using bisection to obtain the fine mesh for cases of $\epsilon=1,0.5,0.1,0.05$, and refine the initial mesh 24 times to obtain the fine mesh for cases of $\epsilon=0.01,0.005$. We set $w=1.0$, $\delta t = 5\times 10^{-3}$ and choose the piecewise linear function as the finite element basis. For the POD type methods, we choose the same parameters as those for solving (\ref{Kol_eq}) expect the coarse meshes used by the two adaptive POD methods. We refine the initial mesh 14 times to obtain the coarse mesh for cases of $\epsilon = 1, 0.5, 0.1$, and refine the initial mesh 15 times to obtain the coarse mesh for cases of $\epsilon = 0.05, 0.01, 0.005$. The time steps corresponding to the coarse finite element spaces are 0.2 and 0.125, respectively. We use the similar number of processors as those for the Kolmogorov flow to do the simulation, that is, we use 72 processors and 180 processors for cases of the fine mesh being obtained by refining the initial mesh 23 times and 24 times, respectively.

\begin{figure}[H]
    \centering
    \subfloat[$\epsilon=1.$]{
        \includegraphics[width=5.5cm]{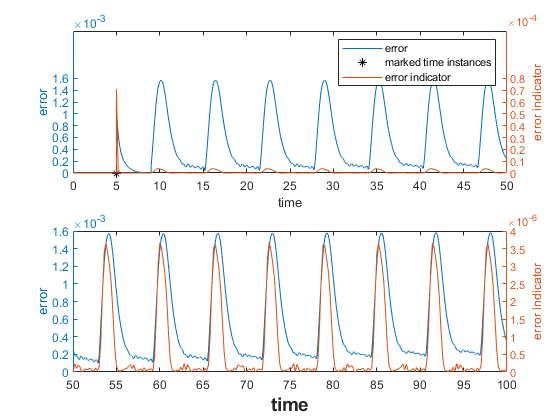}
    }
    \quad
    \subfloat[$\epsilon=0.5.$]{
        \includegraphics[width=5.5cm]{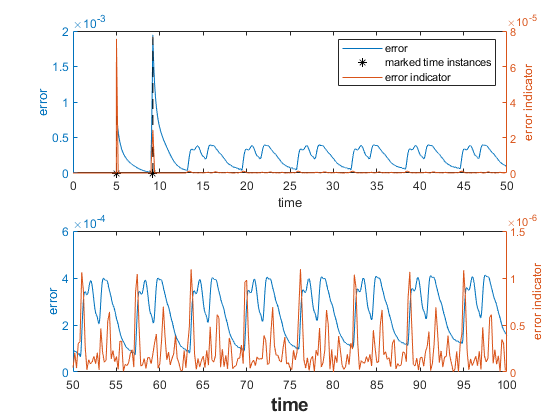}
    }
    \quad
    \subfloat[$\epsilon=0.1.$]{
        \includegraphics[width=5.5cm]{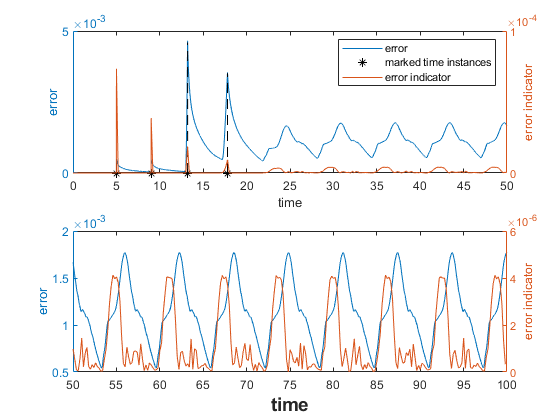}
    }
    \quad
    \subfloat[$\epsilon=0.05.$]{
        \includegraphics[width=5.5cm]{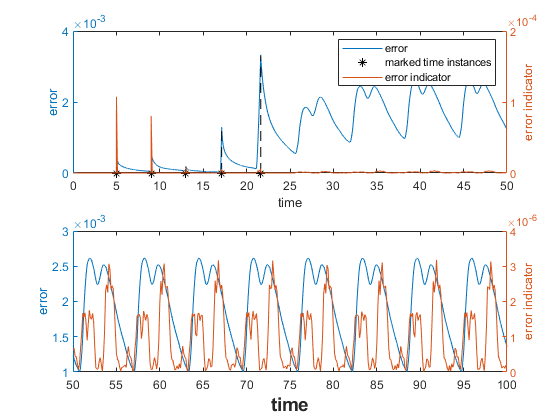}
    }
    \quad
    \subfloat[$\epsilon=0.01.$]{
        \includegraphics[width=5.5cm]{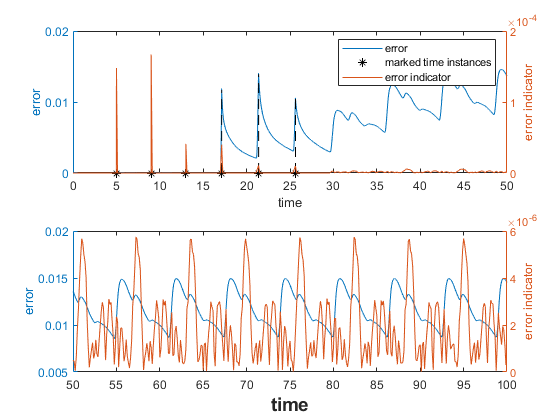}
    }
    \quad
    \subfloat[$\epsilon=0.005.$]{
        \includegraphics[width=5.5cm]{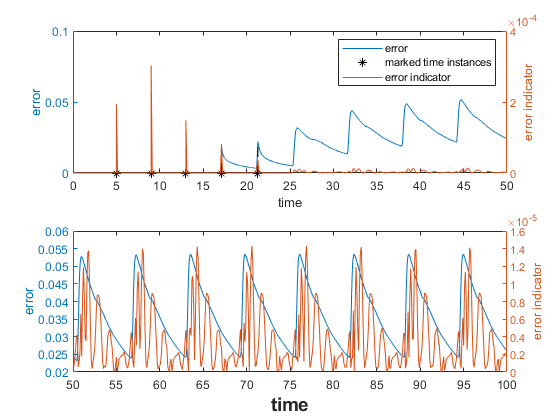}
    }
    \caption{The evolution curves of the error indicator and the error obtined by the method Aug-APOD for solution of (\ref{ABCeq}) with $\epsilon=1,0.5,0.1,0.05,0.01,0.005$, respectively}
    \label{figABC}
\end{figure}

Similarly, we show the error and the error indicator obtained by our method Aug-APOD in Figure \ref{figABC}. The x-axis and y-axis of each sub-figure have the same meaning as those in Figure \ref{figKolmo}. In each sub-figure, the blue curve and orange curve also describe the variation of the error and error indicator, and the black star also denotes the marked time instance. Figure \ref{figABC} shows that the variation of the error indicator is proportional to that of the error. This means that the error indicator is also very effective for this example.
\graphicspath{{Img/Err&ErrInd/}}
The numerical results obtained by methods FEM, POD, TG-APOD and Aug-APOD are shown in Table \ref{tableABC}. The parameters set for the methods TG-APOD and Aug-APOD are the same expect the threshold value $\eta_0$. The notation in Table \ref{tableABC} has the same meaning as those in Table \ref{tableKolmo}.

Similar to the first example, from Table \ref{tableABC}, we see that the POD type methods save much CPU time. The results obtained by the standard POD method are also unreliable as the decrease of $\epsilon$. The two adaptive POD methods are also of high accuracy even for case of $\epsilon=0.005$. Moreover, the adaptive POD methods can save approximately two-thirds of the CPU time compared with that cost by the standard finite element method. This shows that the two adaptive POD methods behave much better than the standard POD method. For the two adaptive POD methods Aug-APOD and TG-APOD, our new method Aug-APOD behaves better than the method TG-APOD. Specifically, the method Aug-APOD can achieve higher accuracy with similar number of POD modes for cases of $\epsilon=0.05,0.01,0.005$. Moreover, for case of $\epsilon=0.005$, the method Aug-APOD can use less CPU time to obtain results as accurate as those obtained by the method TG-APOD.

\begin{figure}[H]
    \centering
    \subfloat[$\epsilon=1.$]{
        \includegraphics[width=5.5cm]{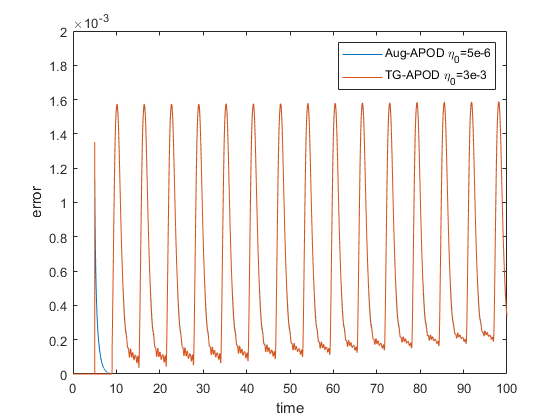}
    }
    \quad
    \subfloat[$\epsilon=0.5.$]{
        \includegraphics[width=5.5cm]{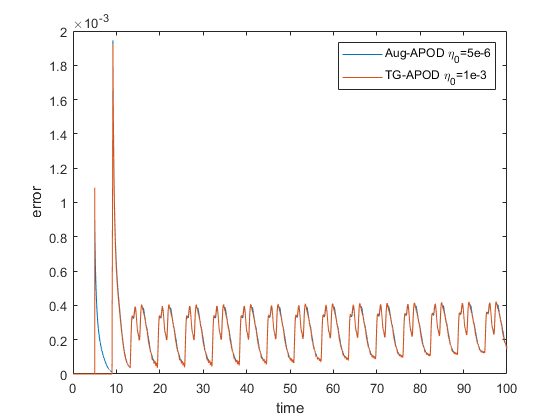}
    }
    \quad
    \subfloat[$\epsilon=0.1.$]{
        \includegraphics[width=5.5cm]{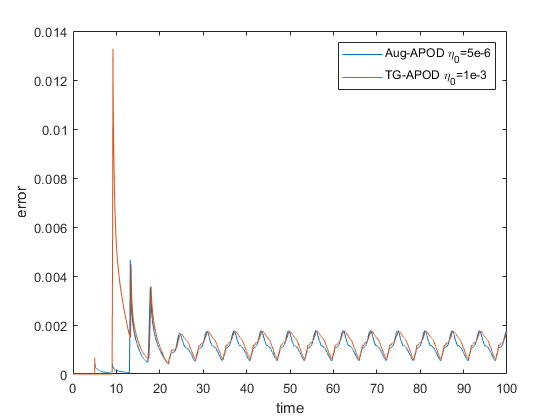}
    }
    \quad
    \subfloat[$\epsilon=0.05.$]{
        \includegraphics[width=5.5cm]{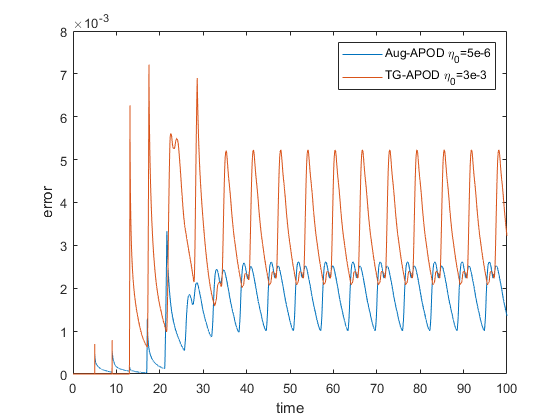}
    }
    \quad
    \subfloat[$\epsilon=0.01.$]{
        \includegraphics[width=5.5cm]{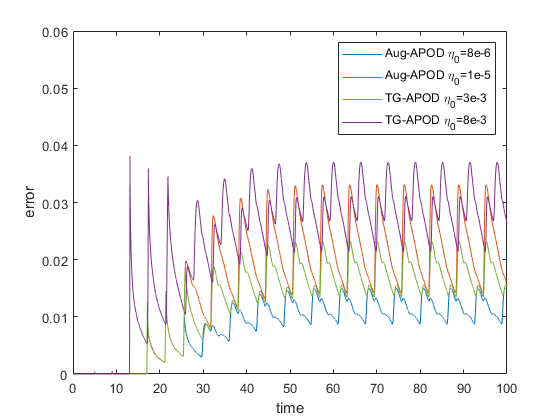}
    }
    \quad
    \subfloat[$\epsilon=0.005.$]{
        \includegraphics[width=5.5cm]{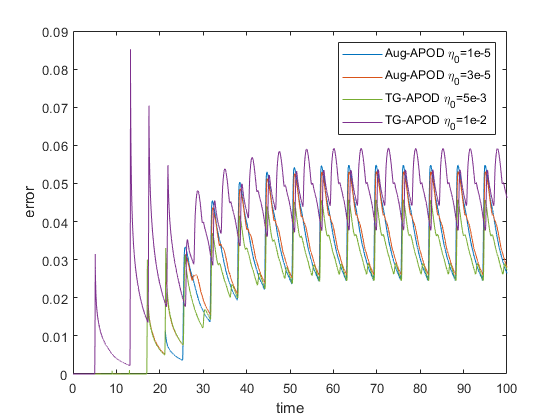}
    }
    \caption{The evolution curves of the error for solution of (\ref{ABCeq}) with different $\epsilon$ by methods TG-APOD and Aug-APOD, respectively}
    \label{figABCTGCos}
\end{figure}

Similarly, more numerical results obtained by using other $\eta_0$ and different coarse meshes for case of $\epsilon= 0.005$ can be found in ``Appendix A ''.
\begin{center}
    \begin{table}[H]
        \caption{The results of (\ref{ABCeq}) with different $\epsilon$ obtained by methods FEM, POD, TG-APOD and Aug-APOD, respectively}\centering
        \label{tableABC}
        \small{
            \begin{tabular}{|cccccccc|}
                \hline
                $\epsilon$ & Method                           & $\eta_{0}$ & Update Times & DOFs     & Error    & Average Error & Time(s)  \\
                \hline
                1          & \textcolor[rgb]{0,0,0}{FEM}      & -          & -            & 10485760 & -        & -             & 20427.20 \\
                           & \textcolor[rgb]{0,0,0}{POD}      & -          & -            & 24       & 0.000974 & 0.013559      & 1681.67  \\
                           & \textcolor[rgb]{0,0,0}{TG-APOD}  & 3e-3       & 1           & 17      & 0.000335 & 0.000517      & 2304.13  \\
                           & \textcolor[rgb]{0,0,0}{Aug-APOD} & 5e-6       & 1           & 17       & 0.000341 & 0.000522      & 2340.55  \\
                           \hline
                0.5        & \textcolor[rgb]{0,0,0}{FEM}      & -          & -            & 10485760 & -        & -             & 18432.76 \\
                           & \textcolor[rgb]{0,0,0}{POD}      & -          & -            & 26       & 0.007569 & 0.069137      & 1625.80  \\
                           & \textcolor[rgb]{0,0,0}{TG-APOD}  & 1e-3       & 2            & 28       & 0.000145 & 0.000234      & 3044.43  \\
                           & \textcolor[rgb]{0,0,0}{Aug-APOD} & 5e-6       & 2            & 28       & 0.000154 & 0.000244      & 3064.83  \\
                           \hline
                0.1        & \textcolor[rgb]{0,0,0}{FEM}      & -          & -            & 10485760 & -        & -             & 17362.39 \\
                           & \textcolor[rgb]{0,0,0}{POD}      & -          & -            & 35       & 0.136214 & 0.373726      & 1727.87  \\
                           & \textcolor[rgb]{0,0,0}{TG-APOD}  & 1e-3       & 4            & 68       & 0.001716 & 0.001258      & 4940.78  \\
                           & \textcolor[rgb]{0,0,0}{Aug-APOD} & 5e-6       & 4            & 67       & 0.001767 & 0.000998      & 4889.35  \\
                           \hline
                0.05       & \textcolor[rgb]{0,0,0}{FEM}      & -          & -            & 10485760 & -        & -             & 17391.44 \\
                           & \textcolor[rgb]{0,0,0}{POD}      & -          & -            & 40       & 0.257052 & 0.514289      & 1814.08  \\
                           & \textcolor[rgb]{0,0,0}{TG-APOD}  & 3e-3       & 5            & 100      & 0.003220 & 0.002760      & 6209.27  \\
                           & \textcolor[rgb]{0,0,0}{Aug-APOD} & 5e-6       & 5            & 99       & 0.001352 & 0.001507      & 6222.50  \\
                           \hline
                0.01       & \textcolor[rgb]{0,0,0}{FEM}      & -          & -            & 16777216 & -        & -             & 16065.81 \\
                           & \textcolor[rgb]{0,0,0}{POD}      & -          & -            & 44       & 0.577333 & 0.697981      & 1515.37  \\
                           & \textcolor[rgb]{0,0,0}{TG-APOD}  & 8e-3       & 5            & 166      & 0.026295 & 0.022662      & 5796.19  \\
                           & \textcolor[rgb]{0,0,0}{TG-APOD}  & 3e-3       & 6            & 192      & 0.015234 & 0.012518      & 6640.86  \\
                           & \textcolor[rgb]{0,0,0}{Aug-APOD} & 1e-5       & 5            & 164      & 0.015639 & 0.016240      & 6037.70  \\
                           & \textcolor[rgb]{0,0,0}{Aug-APOD} & 8e-6       & 6            & 193      & 0.014533 & 0.008428      & 6808.99  \\
                           \hline
                0.005      & \textcolor[rgb]{0,0,0}{FEM}      & -          & -            & 16777216 & -        & -             & 16836.62 \\
                           & \textcolor[rgb]{0,0,0}{POD}      & -          & -            & 44       & 0.812194 & 0.836567      & 1554.84  \\
                           & \textcolor[rgb]{0,0,0}{TG-APOD}  & 1e-2       & 5            & 194      & 0.046166 & 0.039155      & 6437.76  \\
                           & \textcolor[rgb]{0,0,0}{TG-APOD}  & 5e-3       & 6            & 226      & 0.027286 & 0.023730      & 7147.38  \\
                           & \textcolor[rgb]{0,0,0}{Aug-APOD} & 3e-5       & 5            & 193      & 0.027987 & 0.027651      & 6592.85  \\
                           & \textcolor[rgb]{0,0,0}{Aug-APOD} & 1e-5       & 5            & 192      & 0.026378 & 0.027074      & 6517.28  \\
                \hline
            \end{tabular}
        }
    \end{table}
\end{center}
Similar to the first example, we show the comparison between the methods TG-APOD and Aug-APOD in Figure \ref{figABCTGCos}. The x-axis of each figure is the time, the y-axis is the relative error of numerical solution obtained by the adaptive POD methods. 
It is obvious that the error obtained by the method Aug-APOD is smaller than that obtained by the method TG-APOD over the entire time interval for cases of $\epsilon=0.05$ and $\epsilon=0.01$. From these comparison, we can also see that our new method Aug-APOD behaves better than the method TG-APOD. 
\graphicspath{{Img/TGA&Aug_Err/}}
From the two examples we have tested, we see that our new method Aug-APOD has better performance than the exiting methods. 
\section{Concluding remarks}
In this paper, we have proposed a new approach for adaptive POD methods to solve time dependent partial differential equations efficiently and accurately. By augmenting the POD subspace with some auxiliary modes, we have obtained an augmented subspace. We have used the difference between the approximation obtained in this augmented subspace and that obtained in the original POD subspace to construct an error indicator. Using this idea, we have obtained a general framework for the augmented subspace based adaptive POD method. We have then provided two specific methods to obtain the auxiliary mode, one is generating vector randomly, the other is computing the coarse-grid approximations.

We have compared the new approach with the two-grid based adaptive method to test two typical 3D advection-diffusion equations, including both cases with the advection being the Kolmogorov flow and the ABC flow. Numerical results show that our method is more efficient than the existing adaptive methods, especially for cases of $\epsilon=0.01,0.005$. In our future work, we will pay more effort on
applying the method to other types of time dependent partial differential equations and constructing other effective error indicators.

\section*{Acknowledgements}
    This work was supported by the National Key R \& D Program of China under grants 2019YFA0709600 and 2019YFA0709601, the National Natural Science Foundation of China under grants 92270206 and 11671389, and the NSF grants DMS-1952644, DMS-2151235.

\appendix
\section{Numerical experiments for different coarse meshes and different $\eta_0$}
In this section, we use more numerical experiments to illustrate that our new method performs better than the two-grid based adaptive POD method, especially for cases of  $\epsilon=0.01, 0.005$. We denote the two-grid based adaptive POD method and our augmented subspace based adaptive POD method with the coarse grid obtained by refining the initial mesh $N$ times as TG-APOD-$N$ and Aug-APOD-$N$, respectively. The time steps corresponding to different coarse finite spaces are listed in Table \ref{tablegrid}. With a fixed fine grid, we have tested (\ref{Kol_eq}) with $\epsilon=0.01,0.005$ and (\ref{ABCeq}) with $\epsilon = 0.005$ by different coarse meshes and different threshold $\eta_0$, respectively. 
\begin{center}
    \begin{table}[H]
        \caption{\small Time steps corresponding to different coarse meshes}\centering\label{tablegrid}
        \footnotesize{
            \begin{tabular}{|ccccc|}
                \hline
                $N$-Refine  & 14  & 15    & 16  & 17   \\
                \hline
                Time step & 0.2 & 0.125 & 0.1 & 0.05 \\
                \hline
            \end{tabular}
        }
    \end{table}
\end{center}
\subsection{Kolmogorov flow with $\epsilon=0.01$}
The results of (\ref{Kol_eq}) with $\epsilon=0.01$ obtained by different methods and different parameters are shown in Table \ref{tableKolmo_0.01}. From Table \ref{tableKolmo_0.01}, we can draw the similar conclusion as those from Table \ref{tableKolmo}. We first see that all adaptive POD methods can obtain results with some accuracy when compared with those obtained by the FEM method but cost much less CPU time than the FEM method. Then, we see that all adaptive POD methods behave much better than the traditional POD method. With similar number of POD modes, the results obtained by our augmented subspace based adaptive POD method are more accurate than those obtained by other POD type methods. 
\begin{center}
    \begin{table}[H]
        \caption{The results of (\ref{Kol_eq}) with $\epsilon=0.01$ obtained by different POD type methods with different coarse meshes and $\eta_0$, respectively}\centering
        \label{tableKolmo_0.01}
        \small{
            \begin{tabular}{|cccccccc|}
                \hline
                $\epsilon$ & Method                              & $\eta_{0}$ & Update Times & DOFs     & Error    & Average Error & Time(s)  \\
                \hline
                0.01       & \textcolor[rgb]{0,0,0}{FEM}         & -          & -            & 16777216 & -        & -             & 11784.99 \\
                           & \textcolor[rgb]{0,0,0}{POD}         & -          & -            & 43       & 0.919008 & 0.777626      & 1267.46  \\
                           & \textcolor[rgb]{0,0,0}{TG-APOD-14}  & 5e-3       & 5            & 143      & 0.024535 & 0.029794      & 4424.10  \\
                           & \textcolor[rgb]{0,0,0}{TG-APOD-14}  & 3e-3       & 5            & 144      & 0.013786 & 0.017077      & 4522.11  \\
                           & \textcolor[rgb]{0,0,0}{TG-APOD-14}  & 1e-3       & 6            & 167      & 0.009079 & 0.011065      & 5180.53  \\
                           & \textcolor[rgb]{0,0,0}{TG-APOD-15}  & 5e-3       & 5            & 144      & 0.013528 & 0.015620      & 4511.11  \\
                           & \textcolor[rgb]{0,0,0}{TG-APOD-15}  & 3e-3       & 6            & 165      & 0.010431 & 0.010523      & 5043.79  \\
                           & \textcolor[rgb]{0,0,0}{TG-APOD-15}  & 1e-3       & 6            & 166      & 0.007042 & 0.004837      & 5072.23  \\
                           & \textcolor[rgb]{0,0,0}{TG-APOD-16}  & 5e-3       & 5            & 147      & 0.011280 & 0.012568      & 4602.02  \\
                           & \textcolor[rgb]{0,0,0}{TG-APOD-16}  & 3e-3       & 6            & 170      & 0.007904 & 0.008971      & 5130.71  \\
                           & \textcolor[rgb]{0,0,0}{TG-APOD-16}  & 1e-3       & 7            & 196      & 0.005060 & 0.003834      & 5789.31  \\
                \hline
                           & \textcolor[rgb]{0,0,0}{Aug-APOD-14} & 5e-5       & 5            & 149      & 0.022220 & 0.015238      & 4453.28  \\
                           & \textcolor[rgb]{0,0,0}{Aug-APOD-14} & 3e-5       & 5            & 148      & 0.019160 & 0.016858      & 4410.95  \\
                           & \textcolor[rgb]{0,0,0}{Aug-APOD-14} & 1e-5       & 6            & 173      & 0.008892 & 0.004957      & 5302.05  \\
                           & \textcolor[rgb]{0,0,0}{Aug-APOD-14} & 5e-6       & 7            & 199      & 0.001541 & 0.001740      & 5708.77  \\
                           & \textcolor[rgb]{0,0,0}{Aug-APOD-15} & 5e-5       & 5            & 147      & 0.009849 & 0.009620      & 4560.84  \\
                           & \textcolor[rgb]{0,0,0}{Aug-APOD-15} & 3e-5       & 5            & 147      & 0.009588 & 0.009403      & 4505.68  \\
                           & \textcolor[rgb]{0,0,0}{Aug-APOD-15} & 1e-5       & 6            & 172      & 0.006208 & 0.003351      & 5280.25  \\
                           & \textcolor[rgb]{0,0,0}{Aug-APOD-15} & 5e-6       & 7            & 199      & 0.003427 & 0.001915      & 5996.87  \\
                           & \textcolor[rgb]{0,0,0}{Aug-APOD-16} & 5e-5       & 5            & 147      & 0.008976 & 0.009054      & 4557.52  \\
                           & \textcolor[rgb]{0,0,0}{Aug-APOD-16} & 3e-5       & 5            & 147      & 0.008976 & 0.009054      & 4505.68  \\
                           & \textcolor[rgb]{0,0,0}{Aug-APOD-16} & 1e-5       & 6            & 172      & 0.006253 & 0.003377      & 5335.85  \\
                \hline
            \end{tabular}
        }

    \end{table}
\end{center}
\subsection{Kolmogorov flow with $\epsilon=0.005$}
We show the results of (\ref{Kol_eq}) with $\epsilon=0.005$ obtained by different methods and different parameters in Table \ref{tableKolmo_0.005}. For this case of $\epsilon$, we can also draw the similar conclusion as those from Table \ref{tableKolmo} and Table \ref{tableKolmo_0.01}, that is, our augmented subspace based adaptive POD method behave better than other methods.
\begin{center}
    \begin{table}[H]
        \caption{The results of (\ref{Kol_eq}) with $\epsilon=0.005$ obtained by different POD type methods with different coarse meshes and $\eta_0$, respectively}\centering
        \label{tableKolmo_0.005}
        \small{
            \begin{tabular}{|cccccccc|}
                \hline
                $\epsilon$ & Method                              & $\eta_{0}$ & Update Times & DOFs     & Error    & Average Error & Time(s)  \\
                \hline
                0.005      & \textcolor[rgb]{0,0,0}{FEM}         & -          & -            & 16777216 & -        & -             & 12526.50 \\
                           & \textcolor[rgb]{0,0,0}{POD}         & -          & -            & 45       & 1.010868 & 0.871304      & 1306.83  \\
                           & \textcolor[rgb]{0,0,0}{TG-APOD-14}  & 1e-2       & 4            & 139      & 0.110962 & 0.091636      & 4132.52  \\
                           & \textcolor[rgb]{0,0,0}{TG-APOD-14}  & 5e-3       & 5            & 171      & 0.062622 & 0.073958      & 4950.07  \\
                           & \textcolor[rgb]{0,0,0}{TG-APOD-15}  & 1e-2       & 5            & 170      & 0.059025 & 0.037340      & 5040.39  \\
                           & \textcolor[rgb]{0,0,0}{TG-APOD-15}  & 5e-3       & 6            & 200      & 0.018382 & 0.022522      & 5713.98  \\
                           & \textcolor[rgb]{0,0,0}{TG-APOD-16}  & 1e-2       & 5            & 178      & 0.023319 & 0.022661      & 5181.44  \\
                           & \textcolor[rgb]{0,0,0}{TG-APOD-16}  & 5e-3       & 6            & 208      & 0.015640 & 0.014779      & 5887.37  \\
                \hline
                           & \textcolor[rgb]{0,0,0}{Aug-APOD-14} & 5e-5       & 5            & 176      & 0.018007 & 0.018023      & 4893.47  \\
                           & \textcolor[rgb]{0,0,0}{Aug-APOD-14} & 3e-5       & 6            & 206      & 0.013928 & 0.008509      & 5789.30  \\
                           & \textcolor[rgb]{0,0,0}{Aug-APOD-14} & 1e-5       & 6            & 207      & 0.014098 & 0.007749      & 5696.05  \\
                           & \textcolor[rgb]{0,0,0}{Aug-APOD-15} & 5e-5       & 5            & 176      & 0.017100 & 0.017024      & 5084.80  \\
                           & \textcolor[rgb]{0,0,0}{Aug-APOD-15} & 3e-5       & 6            & 207      & 0.013551 & 0.008223      & 5913.76  \\
                           & \textcolor[rgb]{0,0,0}{Aug-APOD-15} & 1e-5       & 7            & 239      & 0.008122 & 0.004791      & 6891.73  \\
                           & \textcolor[rgb]{0,0,0}{Aug-APOD-16} & 5e-5       & 5            & 176      & 0.017162 & 0.016912      & 5295.97  \\
                           & \textcolor[rgb]{0,0,0}{Aug-APOD-16} & 3e-5       & 5            & 176      & 0.017162 & 0.016912      & 5149.72  \\
                           & \textcolor[rgb]{0,0,0}{Aug-APOD-16} & 1e-5       & 7            & 239      & 0.007957 & 0.004724      & 6778.84  \\
                \hline
            \end{tabular}
        }

    \end{table}
\end{center}
\subsection{ABC flow with $\epsilon=0.005$}
In Table \ref{tableABC_0.005}, we show the results of (\ref{ABCeq}) with $\epsilon=0.005$ obtained by different methods and different parameters. Similar to the example of Kolmogorov flow, the adaptive POD methods behave better than other methods, and our augmented subspace based adaptive POD method behave better than the two-grid based adaptive POD method. 
\begin{center}
    \begin{table}[H]
        \caption{The results of (\ref{ABCeq}) with $\epsilon=0.005$ obtained by different POD type methods with different coarse meshes and $\eta_0$, respectively}\centering
        \label{tableABC_0.005}
        \small{
            \begin{tabular}{|cccccccc|}
                \hline
                $\epsilon$ & Method                              & $\eta_{0}$ & Update Times & DOFs     & Error    & Average Error & Times    \\
                \hline
                0.005      & \textcolor[rgb]{0,0,0}{FEM}         & -          & -            & 16777216 & -        & -             & 16836.62 \\
                           & \textcolor[rgb]{0,0,0}{POD}         & -          & -            & 44       & 0.812194 & 0.836567      & 1554.84  \\
                           & \textcolor[rgb]{0,0,0}{TG-APOD-15}  & 1e-2       & 5            & 194      & 0.046166 & 0.039155      & 6437.76  \\
                           & \textcolor[rgb]{0,0,0}{TG-APOD-15}  & 5e-3       & 6            & 226      & 0.027286 & 0.023730      & 7147.38  \\
                           & \textcolor[rgb]{0,0,0}{TG-APOD-16}  & 2e-2       & 5            & 198      & 0.086156 & 0.076253      & 6563.74  \\
                           & \textcolor[rgb]{0,0,0}{TG-APOD-16}  & 1e-2       & 6            & 228      & 0.034736 & 0.035273      & 7395.12  \\
                           & \textcolor[rgb]{0,0,0}{TG-APOD-16}  & 5e-3       & 7            & 259      & 0.025176 & 0.020051      & 8409.76  \\
                           & \textcolor[rgb]{0,0,0}{TG-APOD-17}  & 2e-2       & 5            & 195      & 0.046125 & 0.038615      & 6460.26  \\
                           & \textcolor[rgb]{0,0,0}{TG-APOD-17}  & 1e-2       & 8            & 294      & 0.025939 & 0.019266      & 9540.73  \\
                           & \textcolor[rgb]{0,0,0}{TG-APOD-17}  & 5e-3       & 9            & 323      & 0.017107 & 0.010361      & 10735.34 \\
                \hline

                           & \textcolor[rgb]{0,0,0}{Aug-APOD-15} & 5e-5       & 4            & 159      & 0.065712 & 0.051971      & 5581.35  \\
                           & \textcolor[rgb]{0,0,0}{Aug-APOD-15} & 3e-5       & 5            & 193      & 0.027987 & 0.027651      & 6592.85  \\
                           & \textcolor[rgb]{0,0,0}{Aug-APOD-15} & 2e-5       & 5            & 193      & 0.027987 & 0.027651      & 6480.75  \\
                           & \textcolor[rgb]{0,0,0}{Aug-APOD-15} & 1e-5       & 5            & 192      & 0.026378 & 0.027074      & 6517.28  \\
                           & \textcolor[rgb]{0,0,0}{Aug-APOD-16} & 5e-5       & 5            & 192      & 0.032155 & 0.031447      & 6468.33  \\
                           & \textcolor[rgb]{0,0,0}{Aug-APOD-16} & 3e-5       & 5            & 192      & 0.025921 & 0.026135      & 6512.32  \\
                           & \textcolor[rgb]{0,0,0}{Aug-APOD-16} & 2e-5       & 5            & 192      & 0.025921 & 0.026135      & 6479.52  \\
                           & \textcolor[rgb]{0,0,0}{Aug-APOD-16} & 1e-5       & 9            & 323      & 0.011047 & 0.014037      & 10392.84 \\
                \hline
            \end{tabular}
        }
    \end{table}
\end{center}

\bibliography{bib}

\end{document}